\documentclass[a4paper,11pt,fleqn]{article}
\usepackage{amsmath,amssymb,amsthm,graphicx,subfigure,float,caption,epstopdf,tabularx,color, bm,amsfonts,epic}
\usepackage[top=1in, bottom=1in, left=1.25in, right=1.25in]{geometry}
\usepackage{appendix}
\usepackage{multirow}
\usepackage{diagbox}
\allowdisplaybreaks
\newtheorem{thm}{Theorem}[section]
\newtheorem{lem}{Lemma}[section]
\newtheorem{coro}{Corollary}[section]
\newtheorem{den}{Definition}[section]

\newtheorem{rmk}{Remark}[section]

\newtheorem*{prf}{Proof}
\numberwithin{equation}{section}

\usepackage{indentfirst}
\graphicspath{{Fig}}

\begin{document}
\title{A linearly implicit energy-preserving exponential integrator for the nonlinear Klein-Gordon equation}
\author{Chaolong Jiang$^1$,\ Yushun Wang$^2$ and Wenjun Cai$^2$\footnote{Correspondence author. Email:
caiwenjun@njnu.edu.cn.}\\
{\small $^1$ School of Statistics and Mathematics, }\\
{\small Yunnan University of Finance and Economics, Kunming 650221, P.R. China}\\
{\small $^2$ Jiangsu Key Laboratory for Numerical Simulation of Large Scale Complex Systems,}\\
{\small School of Mathematical Sciences,  Nanjing Normal University,}\\
{\small  Nanjing 210023, P.R. China}\\
}
\date{}
\maketitle

\begin{abstract}
In this paper, we generalize the exponential energy-preserving integrator proposed in the recent paper [SIAM J. Sci. Comput. 38(2016) A1876-A1895] for
conservative systems, which now becomes linearly implicit by further utilizing the idea of the scalar auxiliary variable
approach. Comparing with the original exponential energy-preserving integrator which usually leads to a
nonlinear algebraic system, our new method only involve a linear system with constant coefficient matrix. Taking
the nonlinear Klein-Gordon equation for example, we derive the concrete energy-preserving scheme and
demonstrate its high efficiency through numerical experiments.

\end{abstract}

\section{Introduction}

It is well-known that exponential integrators permit larger step sizes and achieve higher
accuracy than nonexponential ones when the considered problem is a very stiff differential equation such as
highly oscillatory ODEs or semidiscrete time-dependent PDEs. As to exponential integrators, the
earlier attempts can date back to the original paper
by Hersch \cite{Hersch58}, whereas the term ``exponential integrators" was coined in the seminal paper by Hochbruck, Lubich, and Selhofer \cite{HLS98}. Readers are referred to Ref. \cite{HO10an} for details about exponential integrators. Over the years, there has been growing interest in structure-preserving exponential methods, which can preserve as much as possible the physical/geromeric properties of the dynamic system under consideration \cite{ELW06}. Due to the superior properties in the capability for the long-term computation, symplectic exponential methods have attracted much attention (e.g., see Refs. \cite{MW17jcp,SL19jcp,WW18book} and references therein). On the other hand, the energy conservation law is an important property of conservative systems and whether or not can preserve the energy conservation law of the original systems is a criterion to
judge the success of a numerical method for their solution (e.g. see Refs. \cite{FM2011,LQ95,ZGL95} and references therein). Thus, how to design energy-preserving schemes for conservative systems attracts a lot of interest in recent years. The
noticeable ones include the discrete gradient (including averaged vector filed (AVF) method) \cite{CLW17,CGM12,LWQ14,QM08,MQR99}, discrete variational derivative method \cite{DO11,FM2011}, Hamiltonian Boundary
Value Methods (HBVMs) \cite{BIT10}, energy-preserving continuous stage Runge-Kutta (CSRK) methods \cite{H10,MB16,TS12} and local energy-preserving methods \cite{CSjcp2020,GCW14b,JCW19jsc,WWQ08}, and so on. However, to our best knowledge, most existing works on exponential integrators up to now focus on the construction of
explicit schemes and fail to be energy-preserving. In Ref. \cite{CCOfcm08}, Celledoni et al. proposed some implicit and exponential integrators that preserver symmetric and energy of the cubic Schr\"odinger equation by using
the symmetric projection approach \cite{ELW06}. Recently, combining the ideas of exponential integrators and discrete gradients, Li and Wu \cite{LW16b} proposed an energy-preserving exponential scheme for conservative systems, which was revisited and generalized more recently by Shen and Leok \cite{SL19jcp}. Unfortunately,  such scheme is fully implicit. At each time step, one needs to solve a fully nonlinear system and thus it might be very time consuming. Compared with fully implicit schemes, linearly
implicit schemes only require to solve a linear system, which leads to considerably lower
costs than implicit ones \cite{DO11}. As far as we know, there has been no reference considering linearly implicit exponential
schemes for conservative systems, which can inherit the energy.

In this paper, taking the nonlinear Klein-Gordon equation as an example, we propose a novel linearly implicit exponential scheme for conservative systems by combining the ideas of the exponential integrator and the scalar auxiliary variable (SAV) approach \cite{SXY18,SXY17}. The proposed scheme can inherit the energy and enjoy the same computational advantages as the one (see \cite{CJWS19jcp}) provided by the classical SAV approach. The SAV approach  as well as the earlier invariant energy quadratization (IEQ) approach \cite{YZW17,ZYGW17} is developed based on the idea of the energy quadratization, which can result linearly implicit and energy stable schemes for gradient flows. To the best of our  knowledge, there has been no reference considering the combination of the ideas of the SAV approach and the exponential integrator for developing linearly implicit energy-preserving schemes for energy-conserving systems. Taking the nonlinear Klein-Gordon equation for example, we first explore the feasibility.

The outline of this paper is
organized as follows. In Section \ref{Sec:NKGEs:2}, based on idea of the SAV approach, the NKGE \eqref{NKGEs:eq:1.1} is reformulated into an equivalent system  which inherits a modified energy.  In Section \ref{Sec:NKGEs:3}, a second-order centered difference method is applied to the system and we show that the resulting semi-discrete system can preserve the semi-discrete energy. In Section \ref{Sec:NKGEs:4}, a linearly implicit exponential scheme is presented by combining the  exponential integrator and the linearized
Crank-Nicolson method, which inherits the fully discrete energy. Several numerical examples are shown to illustrate
the power of our proposed scheme in Section \ref{Sec:NKGEs:5}. We draw some conclusions in Section \ref{Sec:NKGEs:6}.

\section{Reformation of the model equation through the SAV approach}\label{Sec:NKGEs:2}
The Klein-Gordon equation is frequently used in mathematical models for problems in many fields of science and engineering, particularly in quantum field theory and relativistic quantum mechanics. Here, we consider the following nonlinear Klein-Gordon equation (NKGE)
\begin{align}\label{NKGEs:eq:1.1}
\left\lbrace
  \begin{aligned}
&\partial_{tt}u({\bm x},t)-\omega^2\Delta u({\bm x},t)+G^{'}(u({\bm x},t))=0,\ {\bm x}\in\mathbb{R}^d,\ t>0,\\
&u({\bm x},0)=\phi_1({\bm x}),\ \partial_t u({\bm x},0)=\phi_2({\bm x}),\ {\bm x}\in\mathbb{R}^d,
\end{aligned}\right.
  \end{align}
where $t$ is time variable, ${\bm x}\in\mathbb{R}^d$ is the spatial variable, $u:=u({\bm x},t)$ is a real-valued function, $\omega$ is a real parameter, $\Delta$ is the usual Laplace operator, $G(u)$ is a smooth potential energy function with $G(u)\ge 0$, and $\phi_1:=\phi_1({\bm x})$ and $\phi_2:=\phi_2({\bm x})$ are two given real-value initial data. The NKGE \eqref{NKGEs:eq:1.1} conserves the Hamiltonian energy 
\begin{align}\label{NKGEs:eq:1.2}
E(t)
&=\int_{\mathbb{R}^d}\Big[\frac{1}{2}|\partial_t u|^2+\frac{\omega^2}{2}|\nabla u|^2+G(u)\Big]d{\bm x}=E(0),\ t\ge 0.
\end{align}
In the last few decades, various structure-preserving methods have been developed for solving the NKGE \eqref{NKGEs:eq:1.1}, including symplectic methods ( e.g., see Refs. \cite{Duncan97,McLachlan93,MW17jcp}), multisymplectic methods (e.g., see Refs. \cite{HJLL07,Reich00,SWX12,ZTHTW10}) and energy-preserving methods (e.g., see Refs. \cite{BCI15,CHL08,guo3,WW18iam}),  etc. However, there has been no reference considering a linearly implicit structure-preserving exponential scheme for the NKGE \eqref{NKGEs:eq:1.1} to our knowledge.

 Following the idea of the SAV approach, we introduce a scalar auxiliary variable, as follows:
 \begin{align*}
 q:=q(t)=\sqrt{(G(u),1)+C_0}.
 \end{align*}
 Here $(f,g)$ is the inner product defined by $(f,g)=\int_{\mathbb{R}^d} f\bar{g} d{\bm x}$ where $\bar{g}$ denotes the conjugate of $g$, and $C_0$ is a constant large enough to make $q$ well-posed.
 The Hamiltonian energy \eqref{NKGEs:eq:1.2} is then rewritten as
 \begin{align}\label{NKGEs:eq:2.1}
E(t)=\int_{\mathbb{R}^d}\Big[\frac{1}{2}|\partial_t u|^2+\frac{\omega^2}{2}|\nabla u|^2\Big]d{\bm x}+q^2-C_0.
 \end{align}
 According to the energy variational formula, the NKGE \eqref{NKGEs:eq:1.1} can be reformulated into the following equivalent form
 \begin{align}\label{NKGEs:eq:2.2}
\left\lbrace
  \begin{aligned}
  &\partial_t u = v,\\
  &\partial_t v = \omega^2\Delta u-\frac{G^{'}(u)}{\sqrt{(G(u),1)+C_0}}q,\\
  &\partial_tq = \frac{(G^{'}(u),\partial_t u)}{2\sqrt{(G(u),1)+C_0}},\\
  &u({\bm x},0) = \phi_1({\bm x}),\ \partial_t u({\bm x},0)=\phi_2({\bm x}),\ q(0)=\sqrt{(G(u({\bm x},0)),1)+C_0},
  \end{aligned}\right.\ \ 
  \end{align}
  where ${\bm x}\in\mathbb{R}^d$ and $t>0$.
\begin{thm} The system \eqref{NKGEs:eq:2.2} possesses the following modified energy.
\begin{align}\label{NKGEs:eq:2.3}
E(t)=\int_{\mathbb{R}^d}\Big[\frac{1}{2}|v|^2+\frac{\omega^2}{2}|\nabla u|^2\Big]d{\bm x}+q^2-C_0=E(0),\ t\ge 0.
\end{align}
\end{thm}
   \begin{prf}\rm Taking the inner products with $v$ of the second equality of \eqref{NKGEs:eq:2.2}, we have, together with the first equality of \eqref{NKGEs:eq:2.2}
   \begin{align}\label{NKGEs:eq:2.4}
   \frac{d}{dt}\int_{\mathbb{R}^d}\Big[\frac{1}{2}|v|^2+\frac{\omega^2}{2}|\nabla u|^2\Big]d{\bm x}+\int_{\mathbb{R}^d}\frac{G^{'}(u)\partial_tu}{\sqrt{(G(u),1)+C_0}}q d{\bm x}=0.
   \end{align}
   Multiplying the third equality of \eqref{NKGEs:eq:2.2} by $q$ gives
   \begin{align}\label{NKGEs:eq:2.5}
   \frac{d}{dt}q^2=\int_{\mathbb{R}^d}\frac{G^{'}(u)\partial_t u}{\sqrt{(G(u),1)}}qd{\bm x}.
   \end{align}
   Combining \eqref{NKGEs:eq:2.4} and \eqref{NKGEs:eq:2.5}, one obtains \eqref{NKGEs:eq:2.3} immediately.
    \qed
   \end{prf}

\begin{rmk} The SAV approach can also valid for a more general $G(u)$. Actually, if $G(u)$ is unbounded from below, we can use the
splitting strategy to divide $G(u)$ into several differences which are bounded from below. Then
the energy can be transformed into a quadratic form by introducing multiple scalar auxiliary
variables and the corresponding model reformulation can be obtained (see Ref. \cite{JGCW19}).
\end{rmk}


\section{Energy-preserving spatial semi-discretization}\label{Sec:NKGEs:3}
For simplicity of notation, we shall introduce our scheme in
one space dimension, i.e. $d=1$ in \eqref{NKGEs:eq:2.2}. Generalizations to $d>1$ are straightforward for tensor product grids
and the results remain valid with modifications. For $d=1$, the NKGE \eqref{NKGEs:eq:2.2} is truncated on a bounded interval $(a,b)$ with the periodic boundary condition.

Choose the mesh size $h=(b-a)/N$ with $N$ an even positive integer, and denote the grid points by $x_{j}=jh$ for $j=0,1,2,\cdots,N$; let $u_{j}$ and $v_j$ be the numerical approximations of  $u(x_j,t)$ and $v(x_j,t)$ for $j=0,1,\cdots,N$, respectively, and $u:=(u_{0},u_{1},\cdots,u_{N-1})^{T}$, $v:=(v_{0},v_{1},\cdots,v_{N-1})^{T}$ be the solution vectors and define the following finite difference operators as
\begin{align*}
\delta^{+}_xu_j=\frac{u_{j+1}-u_j}{h},\ \delta_x^2u_j=\frac{u_{j+1}-2u_j+u_{j-1}}{h^2},\ 0\le j\le N-1.
\end{align*}
In addition, for any $u$ and $v$, we define the discrete inner product and notions as follows
\begin{align*}
&\langle u,v\rangle_{l^2}=h\sum_{j=0}^{N-1}u_j\bar{v}_j,\ ||v||_{l^2}^2=h\sum_{j=0}^{N-1}|v_j|^2,\ ||\delta_x^+u||_{l^2}^2=h\sum_{j=0}^{N-1}|\delta_x^+u_j|^2.
\end{align*}

 Then we apply the second-order centered difference scheme for spatial discretization
 \begin{align}\label{NKGEs:eq:3.2}
\left\lbrace
  \begin{aligned}
  &\frac{d}{dt} u=v,\\
  &\frac{d}{dt} v=\omega^2\delta_x^2 u-\frac{G^{'}(u)}{\sqrt{\langle G(u),{\bm 1}\rangle_{l^2}+C_0}}q,\\
  &\frac{d}{dt} q=\frac{\langle G^{'}(u),\frac{d}{dt} u\rangle_{l^2}}{2\sqrt{\langle G(u),{\bm 1}\rangle_{l^2}+C_0}},\\
  &u_j(0)=\phi_1(x_j),\ v_j(0)=\phi_2(x_j),\ q(0)=\sqrt{\langle G(u(0)),{\bm 1}\rangle_{l^2}+C_0},\\
  &u_{0}=u_{N},\  u_{-1}=u_{N-1},
  \end{aligned}\right.\ \ 
  \end{align}
where $G(u)=(G(u_0,u_1,\cdots,u_{N-1}))^T$ and $0\le j\le N$.

 \begin{thm} The semi-discrete system \eqref{NKGEs:eq:3.2} admits the semi-discrete modified energy
\begin{align}\label{NKGEs:eq:3.3}
E_h(t)=\frac{1}{2}||v||_{l^2}^2+\frac{\omega^2}{2}||\delta_x^+u||_{l^2}^2+q^2-C_0=E_h(0),\ t\ge 0.
\end{align}
\end{thm}
   \begin{prf}\rm Taking the discrete inner products with $v$ of the second equality of \eqref{NKGEs:eq:3.2}, we have, together with the first equality of \eqref{NKGEs:eq:3.2}
   \begin{align}\label{NKGEs:eq:3.4}
   \frac{d}{dt}\Big[\frac{1}{2}||v||^2+\frac{\omega^2}{2}||\delta_x^+ u||_{l^2}^2\Big]+\frac{\langle G^{'}(u),\frac{d}{dt}u\rangle_{l^2}}{\sqrt{\langle G(u),{\bm 1}\rangle_{l^2}+C_0}}q=0.
   \end{align}
   Multiplying with the third equality of \eqref{NKGEs:eq:3.2} by $q$  reads
   \begin{align}\label{NKGEs:eq:3.5}
   \frac{d}{dt}q^2=\frac{\langle G^{'}(u),\frac{d}{dt} u\rangle_{l^2}}{\sqrt{\langle G(u),{\bm 1}\rangle_{l^2}+C_0}}q.
   \end{align}
   Combining \eqref{NKGEs:eq:3.4} and \eqref{NKGEs:eq:3.5}, one obtains \eqref{NKGEs:eq:3.3} immediately.
   \qed
   \end{prf}
\section{Construction of the linearly implicit energy-preserving exponential scheme }\label{Sec:NKGEs:4}

Choose $\tau$ be the time step, and denote $t_{n}=n\tau$ for $n=0,1,2\cdots$; let $u_{j}^n$ be the numerical approximation of  $u(x_j,t_n)$ for $j=0,1,\cdots,N$ and $n=0,1,2,\cdots$; denote $u^n$ as the solution
vector at $t=t_n$ and define
\begin{align*}
&\delta_{t} {u}_j^{n}
=\frac{{u}_j^{n+1}-{u}_j^{n}}{\tau},\ {u}_j^{n+\frac{1}{2}}=\frac{{u}_j^{n+1}+{u}_j^{n}}{2},\ \hat{u}_j^{n+\frac{1}{2}}=\frac{3{u}_j^{n}-{u}_j^{n-1}}{2}, 0\le j\le N-1.
\end{align*}
 {\begin{den}\label{ESAV-den-4.1} Throughout this paper, for a given sufficiently smooth function $f$ in the neighborhood of zero ($f(0):= \displaystyle\lim_{x \to 0} f(x)$ when 0 is a removable singularity),
  \begin{align*}
  f(x)=\sum_{k=0}^{\infty}\frac{f^{(k)}(0)}{k!}x^k,
  \end{align*}
  and for a matrix $A$, the matrix-valued function is defined by
  \begin{align*}
  f(A)= \sum_{k=0}^{\infty}\frac{f^{(k)}(0)}{k!}A^k.
  \end{align*}
  For more details about functions of matrices, please refer to Ref. \cite{Higham2008}.
  \end{den}}

Let $z(t)=(u(t),v(t))^T,\ f(u(t),q(t))=(\frac{G^{'}(u(t))}{\sqrt{\langle G(u(t)),{\bm 1}\rangle_{l^2}+C_0}}q(t),0)^T$ and
\begin{align*}
 S=\left(\begin{array}{cc}
              O& I \\
              -I& O\\
             \end{array}
\right),\ M=\left(\begin{array}{cc}
              -\omega^2 B_2& O \\
              O& I\\
             \end{array}
\right).
\end{align*}
Here, matrix $B_2$ represents the operator $\delta_x^2$ under the periodic boundary condition. In addition, it holds \cite{HNO06}
\begin{align}\label{ESAV-B-2-matrix}
B_2=F^H\Lambda F,\ \Lambda=\text{diag}\Big[\lambda_{0},\lambda_{1},\cdots,\lambda_{N-1}\Big],\ \lambda_{j}=-\frac{4}{h^{2}}\sin^{2}\frac{j\pi}{N},
\end{align}
where ${F}$ is the discrete Fourier matrix of order $N$ and ${F}^H$ represents the conjugate transpose of ${F}$.

Integrating the equation \eqref{NKGEs:eq:3.2} from $t_n$ to $t_{n+1}$, we then have
\begin{align}\label{NKGEs:eq:3.6}
&z(t_n+\tau)=\exp(V)z(t_n)+\tau\int_0^1\exp((1-\xi)V)Sf(u(t_n+\xi\tau),q(t_n+\xi\tau))d\xi,\\\label{NKGEs:eq:3.7}
&q(t_n+\tau)=q(t_n)+\tau\int_0^1\frac{\langle G^{'}(u(t_n+\xi\tau)),\frac{d}{dt}u\rangle_{l^2}}{2\sqrt{\langle G(u(t_n+\xi\tau)),{\bm 1}\rangle_{l^2}+C_0}}d\xi,
\end{align}
where $V=\tau SM$.

Replacing $f(u(t_n+\xi\tau),q(t_n+\xi\tau))$ and $\frac{\langle G^{'}(u(t_n+\xi\tau)),\frac{d}{dt}u\rangle_{l^2}}{2\sqrt{\langle G(u(t_n+\xi\tau)),{\bm 1}\rangle_{l^2}+C_0}}$ with the linearized Crank-Nicolson method $f(\hat{u}^{n+\frac{1}{2}},q^{n+\frac{1}{2}})$ and $\frac{\langle G^{'}(\hat{u}^{n+\frac{1}{2}}),\delta_t u^n\rangle_{l^2}}{2\sqrt{\langle G(\hat{u}^{n+\frac{1}{2}}),{\bm 1}\rangle_{l^2}+C_0}}$, respectively, and we obtain the new scheme, as follows: 
\begin{align}\label{NKGEs:eq:3.8}
&z^{n+1}=\exp(V)z^n+\tau\phi(V)Sf(\hat{u}^{n+\frac{1}{2}},q^{n+\frac{1}{2}}),\\\label{NKGEs:eq:3.9}
&q^{n+1}=q^n+\tau\frac{\langle G^{'}(\hat{u}^{n+\frac{1}{2}}),\delta_t u^n\rangle_{l^2}}{2\sqrt{\langle G(\hat{u}^{n+\frac{1}{2}}),{\bm 1}\rangle_{l^2}+C_0}},
\end{align}
 where
\begin{align*}
\phi(V)=\int_0^1\exp((1-\xi)V)d\xi,\ f(\hat{u}^{n+\frac{1}{2}},q^{n+\frac{1}{2}})=\Big(\frac{G^{'}(\hat{u}^{n+\frac{1}{2}})^T}{\sqrt{\langle G(\hat{u}^{n+\frac{1}{2}}),{\bm 1}\rangle_{l^2}+C_0}},0\Big)^T,
\end{align*}
 and $n=1,2,\cdots$. The initial and boundary conditions in \eqref{NKGEs:eq:3.2} are discretized as
\begin{align*}
&u_j^0=\phi_1(x_j),\ v_j^0=\phi_2(x_j),\ q^0=\sqrt{\langle G(u^0),{\bm 1}\rangle_{l^2}+C_0},\ j=0,1,2,\cdots,N,\\
&u_{0}^n=u_{N}^n,\  u_{-1}^n=u_{N-1}^n,\ n\ge 0.
\end{align*}

  \begin{rmk}\label{ESAV-remark-4.2} Since the proposed scheme \eqref{NKGEs:eq:3.8}-\eqref{NKGEs:eq:3.9} is a three level scheme, we calculate ${z}^1$ and $q^1$ by using ${u}^0$ instead of $\hat{ u}^{\frac{1}{2}}$ for the first step.
  \end{rmk}

Then, we show that the scheme \eqref{NKGEs:eq:3.8}-\eqref{NKGEs:eq:3.9} can preserve the fully discrete modified energy. To begin with, we give the following preliminary Lemma presented in Ref. \cite{LW16b}.

\begin{lem}\label{2SG:lem:4.1}
For any symmetric matrix $M$, and scalar $\tau\ge 0$, the matrix $$ A=\exp(V)^{T}M\exp(V)-M$$ is a nilpotent matrix, when
$S$ is skew symmetric.

\end{lem}

We next have the result, as follows:
  \begin{thm}\label{ESAV-KNLS-4.1} The proposed scheme \eqref{NKGEs:eq:3.8}-\eqref{NKGEs:eq:3.9} can preserve the following discrete modified energy
\begin{align}\label{NKGEs:3.12}
E_h^{n+1}=E_h^n,\ E_h^n=\frac{1}{2}||v^n||_{l^2}^2+\frac{\omega^2}{2}||\delta_x^+u^n||_{l^2}^2+(q^n)^2-C_0,
\end{align}
for $n=0,1,2,\cdots.$
\end{thm}
   \begin{prf}\rm We first note that the matrix $M$ is singular, and assume that $\{M_{\epsilon}\}$ is a series of symmetric and nonsingular matrices, which converge to $M$ when $\epsilon\rightarrow 0$. Let $z_{\epsilon}^{n}$ and $q_{\epsilon}^{n}$ satisfy the perturbed scheme
   \begin{align}\label{NKGEs:eq-3.12}
&z_{\epsilon}^{n+1}=\exp(V_{\epsilon})z_{\epsilon}^n+\tau\phi(V_{\epsilon})Sf(\hat{u}_{\epsilon}^{n+\frac{1}{2}},q_{\epsilon}^{n+\frac{1}{2}}),\\\label{NKGEs:eq-3.13}
&q_{\epsilon}^{n+1}=q_{\epsilon}^n+\tau\frac{\langle G^{'}(\hat{u}_{\epsilon}^{n+\frac{1}{2}}),\delta_t u_{\epsilon}^n\rangle_{l^2}}{2\sqrt{\langle G(\hat{u}_{\epsilon}^{n+\frac{1}{2}}),{\bm 1}\rangle_{l^2}+C_0}},
\end{align} where $V_{\epsilon}=\tau SM_{\epsilon}$ and $n=1,2,\cdots$. Denote $\tilde{f}_{\epsilon}:=M_{\epsilon}^{-1}f_{\epsilon}=M_{\epsilon}^{-1}f(\hat{u}_{\epsilon}^{n+\frac{1}{2}},q_{\epsilon}^{n+\frac{1}{2}})$ and
\begin{align}\label{NKGEs-eq-3.14}
{E_{\epsilon,h}^{n}}=\frac{h}{2}(z_{\epsilon}^{n})^TM_{\epsilon}z_{\epsilon}^{n}+ (q_{\epsilon}^{n})^2.
\end{align}
 Then we have
   \begin{align}\label{NKGEs:eq-3.14}
   &\frac{1}{2}(z_{\epsilon}^{n+1})^TM_{\epsilon}z_{\epsilon}^{n+1}\nonumber\\
   &=\frac{1}{2}\Big[(z_{\epsilon}^n)^T\exp(V_{\epsilon})^T+\tau f_{\epsilon}^TS^T\phi(V_{\epsilon})^T\Big]M_{\epsilon}\Big[\exp(V_{\epsilon})z_{\epsilon}^n
 +\tau\phi(V_{\epsilon})Sf_{\epsilon}\Big]\nonumber\\
   &=\frac{1}{2}(z_{\epsilon}^{n})^T\exp(V_{\epsilon})^TM_{\epsilon}\exp(V_{\epsilon})z_{\epsilon}^n+(z_{\epsilon}^{n})^T\exp(V_{\epsilon})^TM_{\epsilon}
   \Big[\exp(V_{\epsilon})-I\Big]\tilde{f}_{\epsilon}\nonumber\\
   &~~~+\frac{1}{2}\tilde{f}_{\epsilon}^T\Big[\exp(V_{\epsilon})^T-I\Big]M_{\epsilon}\Big[\exp(V_{\epsilon})-I\Big]\tilde{f}_{\epsilon}\nonumber\\
   &=\frac{1}{2}(z_{\epsilon}^{n})^T\exp(V_{\epsilon})^TM_{\epsilon}\exp(V_{\epsilon})z_{\epsilon}^n+(z_{\epsilon}^{n})^T\Big[\exp(V_{\epsilon})^TM_{\epsilon}
   \exp(V_{\epsilon})-\exp(V_{\epsilon})^TM_{\epsilon}\Big]\tilde{f}_{\epsilon}\nonumber\\
   &~~~+\frac{1}{2}\tilde{f}_{\epsilon}^T\Big[\exp(V_{\epsilon})^TM_{\epsilon}\exp(V_{\epsilon})-\exp(V_{\epsilon})^TM_{\epsilon}
   -M_{\epsilon}\exp(V_{\epsilon})+M_{\epsilon}\Big]\tilde{f}_{\epsilon}.
   \end{align}
   On the other hand, it follows from \eqref{NKGEs:eq-3.13} that
   \begin{align}\label{NKGEs:eq-3.15}
  (q_{\epsilon}^{n+1})^2-(q_{\epsilon}^n)^2&=\frac{\langle G^{'}(\hat{u}_{\epsilon}^{n+\frac{1}{2}}),u_{\epsilon}^{n+1}-u_{\epsilon}^n\rangle_{l^2}}{\sqrt{\langle G(\hat{u}_{\epsilon}^{n+\frac{1}{2}}),1\rangle_{l^2}+C_0}}q_{\epsilon}^{n+\frac{1}{2}}\nonumber\\
  &=h((z_{\epsilon}^{n+1})^T-(z_{\epsilon}^n)^T)f_{\epsilon}\nonumber\\
  &=h(z_{\epsilon}^n)^T\Big[\exp(V_{\epsilon})^T-I\Big]f_{\epsilon}+\tau hf_{\epsilon}^TS^T\phi(V_{\epsilon})^Tf_{\epsilon}\nonumber\\
  &=h(z_{\epsilon}^n)^T\Big[\exp(V_{\epsilon})^TM_{\epsilon}-M_{\epsilon}\Big]\tilde{f}_{\epsilon}+h\tilde{f}_{\epsilon}^TV_{\epsilon}^T
  \phi(V_{\epsilon})^TM_{\epsilon}\tilde{f}_{\epsilon}\nonumber\\
  &=h(z_{\epsilon}^n)^T\Big[\exp(V_{\epsilon})^TM_{\epsilon}-M_{\epsilon}\Big]\tilde{f}_{\epsilon}+h\tilde{f}_{\epsilon}^T\Big[\exp(V_{\epsilon})^TM_{\epsilon}
  -M_{\epsilon}\Big]\tilde{f}_{\epsilon}.
   \end{align}
   Then, we can deduce from \eqref{NKGEs:eq-3.14} and \eqref{NKGEs:eq-3.15} that
   \begin{align*}
   &E_{\epsilon,h}^{n+1}-E_{\epsilon,h}^{n}\nonumber\\
   &=\frac{h}{2}(z_{\epsilon}^{n+1})^TM_{\epsilon}z_{\epsilon}^{n+1}-\frac{h}{2}(z_{\epsilon}^{n})^TM_{\epsilon}z_{\epsilon}^{n}+(q_{\epsilon}^{n+1})^2-(q_{\epsilon}^n)^2\nonumber\\
   &=\frac{h}{2}(z_{\epsilon}^{n})^T\Big[\exp(V_{\epsilon})^TM_{\epsilon}\exp(V_{\epsilon})-M_{\epsilon}\Big]z_{\epsilon}^n
   +h(z_{\epsilon}^{n})^T\Big[\exp(V_{\epsilon})^TM_{\epsilon}
   \exp(V_{\epsilon})-M_{\epsilon}\Big]\tilde{f}_{\epsilon}\nonumber\\
   &~~~+\frac{h}{2}\tilde{f}_{\epsilon}^T\Big[\exp(V_{\epsilon})^TM_{\epsilon}\exp(V_{\epsilon})-M_{\epsilon}\Big]\tilde{f}_{\epsilon}
   +\frac{h}{2}\tilde{f}_{\epsilon}^T\Big[\exp(V_{\epsilon})^TM_{\epsilon}
   -M_{\epsilon}\exp{V_{\epsilon}}\Big]\tilde{f}_{\epsilon}\nonumber\\
   &=\frac{h}{2}(z_{\epsilon}^{n}+{\tilde{f}_{\epsilon}})^TA_{\epsilon}(z_{\epsilon}^n+{\tilde{f}_{\epsilon}})
   +\frac{h}{2}{\tilde{f}_{\epsilon}^T}C_{\epsilon}{\tilde{f}_{\epsilon}}=0,
   \end{align*}
   where $A_{\epsilon}=\exp(V_{\epsilon})^{T}M_{\epsilon}\exp(V_{\epsilon})-M_{\epsilon}$ and $C_{\epsilon}=\exp(V_{\epsilon})^{T}M_{\epsilon}-M_{\epsilon}\exp(V_{\epsilon})$. The last equality
is from Lemma \ref{2SG:lem:4.1} and the skew symmetry of the matrix $C_{\epsilon}$.
Thus, when $\epsilon\rightarrow 0$, $z_{\epsilon}^n\rightarrow z^n$, $q_{\epsilon}^n\rightarrow q^n$ and \eqref{NKGEs-eq-3.14} lead to
\begin{align*}
E_{h}^{n+1}=E_{h}^{n}.
\end{align*}
   This completes the proof.
   \qed
   \end{prf}
\begin{coro} Supposing $\phi_1\in H^1(\mathbb{R})$ and $\phi_2\in L^2(\mathbb{R})$, it then follows from \eqref{NKGEs:3.12} that
\begin{align*}
||v^n||_{l^2}\le C,\ ||\delta_x^+u^n||_{l^2}\le C,\ |q^n|\le C,\ n=1,2,\cdots,
\end{align*}
which implies that the proposed scheme is unconditionally stable.
\end{coro}
{\begin{rmk} It should be remarked that there have been various works dedicated to deriving energy stable schemes based on exponential time integrations for gradient flows in recent years (e.g., see Refs. \cite{CLLWW19,DJLQsima19,JLQZmc18}). However, such schemes cannot be directly extended to construct energy-preserving schemes for general conservative systems since explicit approximations of the temporal integral of the nonlinear term  do not satisfy a discrete analog of the chain rule that ensures energy conservation.
\end{rmk}}

  Besides its energy-preserving property, a most remarkable thing about the above scheme is that it can be solved efficiently. We rewrite \eqref{NKGEs:eq:3.8} and \eqref{NKGEs:eq:3.9} as
  \begin{align}\label{NKGEs:eq:3.13}
&u^{n+1}=\exp_{11}u^n+\exp_{12}v^n-\tau\phi_{12}\frac{G^{'}(\hat{u}^{n+\frac{1}{2}})}{\sqrt{\langle G(\hat{u}^{n+\frac{1}{2}}),{\bm 1}\rangle_{l^2}+C_0}}{q}^{n+\frac{1}{2}},\\\label{NKGEs:eq:3.14}
&v^{n+1}=\exp_{21}u^n+\exp_{22}v^n-\tau\phi_{22}\frac{G^{'}(\hat{u}^{n+\frac{1}{2}})}{\sqrt{\langle G(\hat{u}^{n+\frac{1}{2}}),{\bm 1}\rangle_{l^2}+C_0}}{q}^{n+\frac{1}{2}},\\\label{NKGEs:eq:3.15}
&q^{n+1}=q^n+\tau\frac{\langle G^{'}(\hat{u}^{n+\frac{1}{2}}),\delta_t u^n\rangle_{l^2}}{2\sqrt{\langle G(\hat{u}^{n+\frac{1}{2}}),{\bm 1}\rangle_{l^2}+C_0}},
\end{align}
where $\exp(V)$ and $\phi(V)$ are partitioned into
\begin{align*}
\exp(V)=\left(\begin{array}{cc}
              \exp_{11}&\exp_{12} \\
              \exp_{21}&\exp_{22}\\
             \end{array}
\right),\
\phi(V)=\left(\begin{array}{cc}
              \phi_{11}& \phi_{12} \\
             \phi_{21}& \phi_{22}\\
             \end{array}
\right).
\end{align*}
Next, by eliminating $q^{n+\frac{1}{2}}$ in \eqref{NKGEs:eq:3.13} , we have
  \begin{align}\label{NKGEs:eq:3.16}
  u^{n+1}+\gamma \langle G^{'}(\hat{u}^{n+\frac{1}{2}}),u^{n+1} \rangle_{l^2}={g}^n,
  \end{align}
  where
  \begin{align*}
  \gamma=\frac{\tau\phi_{12}G^{'}(\hat{u}^{n+\frac{1}{2}})}{4\langle G(\hat{u}^{n+\frac{1}{2}}),{\bm 1}\rangle_{l^2}+4C_0},
  \end{align*} and
  \begin{align*}
   {g}^n=\exp_{11}u^n+\exp_{12}v^n-\tau\phi_{12}\frac{G^{'}(\hat{u}^{n+\frac{1}{2}})}{\sqrt{\langle G(\hat{u}^{n+\frac{1}{2}}),1\rangle_{l^2}+C_0}}{q}^{n}+\gamma\langle G^{'}(\hat{u}^{n+\frac{1}{2}}),u^n\rangle_{l^2}.
  \end{align*}
  We take the discrete inner product of \eqref{NKGEs:eq:3.16} with $G^{'}(\hat{u}^{n+\frac{1}{2}})$ and have
    \begin{align*}
   \Big(1+\langle G^{'}(\hat{u}^{n+\frac{1}{2}}),\gamma\rangle_{l^2}\Big)\langle G^{'}(\hat{u}^{n+\frac{1}{2}}),u^{n+1} \rangle_{l^2}=\langle G^{'}(\hat{u}^{n+\frac{1}{2}}), {g}^n\rangle_{l^2}.
  \end{align*}
Notice $\langle G^{'}(\hat{u}^{n+\frac{1}{2}}),\gamma\rangle_{l^2}\ge 0$, since $\phi_{12}$ is a symmetrical positive semidefinite matrix. We then obtain from the above that
    \begin{align}\label{NKGEs:eq:3.17}
   \langle G^{'}(\hat{u}^{n+\frac{1}{2}}),u^{n+1} \rangle_{l^2}=\frac{\langle G^{'}(\hat{u}^{n+\frac{1}{2}}), {g}^n\rangle_{l^2}}{1+\langle G^{'}(\hat{u}^{n+\frac{1}{2}}),\gamma\rangle_{l^2}}.
  \end{align}
After solving $\langle G^{'}(\hat{u}^{n+\frac{1}{2}}),u^{n+1} \rangle_{l^2}$ from the linear system \eqref{NKGEs:eq:3.17}, ${u}^{n+1}$ is then updated from \eqref{NKGEs:eq:3.16}. Subsequently, $q^{n+1}$ is obtained from \eqref{NKGEs:eq:3.15}. Finally, we get $v^{n+1}$ from \eqref{NKGEs:eq:3.14}.

{\begin{rmk}\label{ESAV-rk-4.2} On the one hand, $\exp(V)$ and $\phi(V)$ can be efficiently implemented via fast Fourier transform since the calculation of exponentials is time-consuming in general. Actually, according to Definition \ref{ESAV-den-4.1} and \eqref{ESAV-remark-4.2}, we have
\begin{align*}
\exp(V)&=I+V+\frac{V^2}{2!}+\cdots+\frac{V^k}{k!}+\cdots\\
&=\Bigg(\begin{matrix}
              F^H\cosh(\tau\omega\Lambda^{\frac{1}{2}})F&\ F^H(\omega\Lambda^{\frac{1}{2}})^{-1}\sinh(\tau\omega\Lambda^{\frac{1}{2}})F \\
              F^H\omega\Lambda^{\frac{1}{2}}\sinh(\tau\omega\Lambda^{\frac{1}{2}})F &F^H\cosh(\tau\omega\Lambda^{\frac{1}{2}})F\\
             \end{matrix}
\Bigg).
\end{align*}
By the similar argument as above, we obtain
\begin{align*}
\phi(V)=\left(\begin{array}{cc}
F^{H}\tau^{-1}(\omega\Lambda^{\frac{1}{2}})^{-1}\sinh(\tau\omega\Lambda^{\frac{1}{2}})F&F^{H}(\tau\omega^2\Lambda)^{-1}(\cosh(\tau\omega\Lambda^{\frac{1}{2}})-I)F\\
F^{H}\tau^{-1}(\cosh(\tau\omega\Lambda^{\frac{1}{2}})-I)F&F^{H}\tau^{-1}(\omega\Lambda^{\frac{1}{2}})^{-1}\sinh(\tau\omega\Lambda^{\frac{1}{2}})F\\
\end{array}
\right).
\end{align*}
\end{rmk}
Here, we should note that
\begin{align*}
&\Lambda^{\frac{1}{2}}=\text{diag}\Big[{\lambda^{\frac{1}{2}}_{0}},{\lambda^{\frac{1}{2}}_{1}},\cdots,{\lambda^{\frac{1}{2}}_{N-1}}\Big],\\ &(\omega\Lambda^{\frac{1}{2}})^{-1}\sinh(\tau\omega\Lambda^{\frac{1}{2}})=\text{\rm diag}\Bigg[\tau,\frac{\sinh(\tau\omega\lambda_1^{\frac{1}{2}})}{\omega\lambda_1^{\frac{1}{2}}},
\cdots,\frac{\sinh(\tau\omega\lambda_{N-1}^{\frac{1}{2}})}{\omega\lambda_{N-1}^{\frac{1}{2}}}\Bigg],\\ &(\tau\omega^2\Lambda)^{-1}(\cosh(\tau\omega\Lambda^{\frac{1}{2}})-I)=\text{\rm diag}\Big[\frac{\tau}{2},\frac{\cosh(\tau\omega\lambda_1^{\frac{1}{2}})-1}{\tau\omega^2\lambda_1}, \cdots,\frac{\cosh(\tau\omega\lambda_{N-1}^{\frac{1}{2}})-1}{\tau\omega^2\lambda_{N-1}}\Big].
\end{align*}
On the other hand, small modifications would allow us to efficiently implement $\exp(V)$ and $\phi(V)$ in two dimensional case.
}

\section{Numerical examples}\label{Sec:NKGEs:5}

In the previous sections, set the nonlinear Klein-Gordon equation as an example, we present a novel linearly implicit energy-preserving exponential integrator
(denoted by ESAVS) for the conservative system. In this section, we report the numerical performance in accuracy, CPU time and
energy preservation of the energy-preserving exponential integrator scheme  for the nonlinear Klein-Gordon equation and the nonlinear Schr\"odinger equation, respectively. Furthermore, we compare the proposed scheme with the exponential averaged vector filed scheme (denoted by EAVFS) proposed in Ref. \cite{LW16b}. All computations are carried out via Matlab 7.0 with AMD A8-7100 and RAM 4GB. In addition, the standard fixed-point iteration is used for EAVFS and the iteration will terminate when the infinity norm of the error between two adjacent iterative steps is less than $10^{-14}$. In order to quantify the numerical solution, we use the $L^2$- and $L^{\infty}$-norms of the error between the numerical solution $u_j^n$ and the exact solution
$u(x_j,t_n)$, respectively, as
\begin{align*}
e_{h,2}(t_n)=\Bigg({h\sum_{j=0}^{N-1}|u_{j}^n-u(x_j,t_n)|^2}\Bigg)^{\frac{1}{2}},\ e_{h,\infty}(t_n)=\max\limits_{0\le j\le N-1}|u_{j}^n-u(x_j,t_n)|,\ n\ge 0.
\end{align*} 

\subsection{Nonlinear Klein-Gordon equation}
 We first consider the one dimensional nonlinear sine-Gordon equation as follows:
\begin{align}\label{NKGEs-4.1}
&\partial_{tt}u-\partial_{xx} u+\sin(u)=0,\ x\in\mathbb{R},\ t>0,
\end{align}
with initial conditions
\begin{align*}
&u(x,0)=0,\ u_t(x,0)=4\text{sech}(x),\ x\in\mathbb{R}.
\end{align*}
Equation \eqref{NKGEs-4.1} possesses the analytical solution
\begin{align*}
u(x,t)=4\arctan(t\text{sech}(x)),\ x\in\mathbb{R},\ t\ge 0.
\end{align*}

First of all, we present the time mesh refinement tests to show the order of accuracy
of the proposed scheme. We choose the parameter $C_0=1$ and the computational domain $\Omega=[-20,20]$ with a periodic boundary condition.

\begin{figure}[H]
\centering
\begin{minipage}[t]{65mm}
\includegraphics[width=65mm]{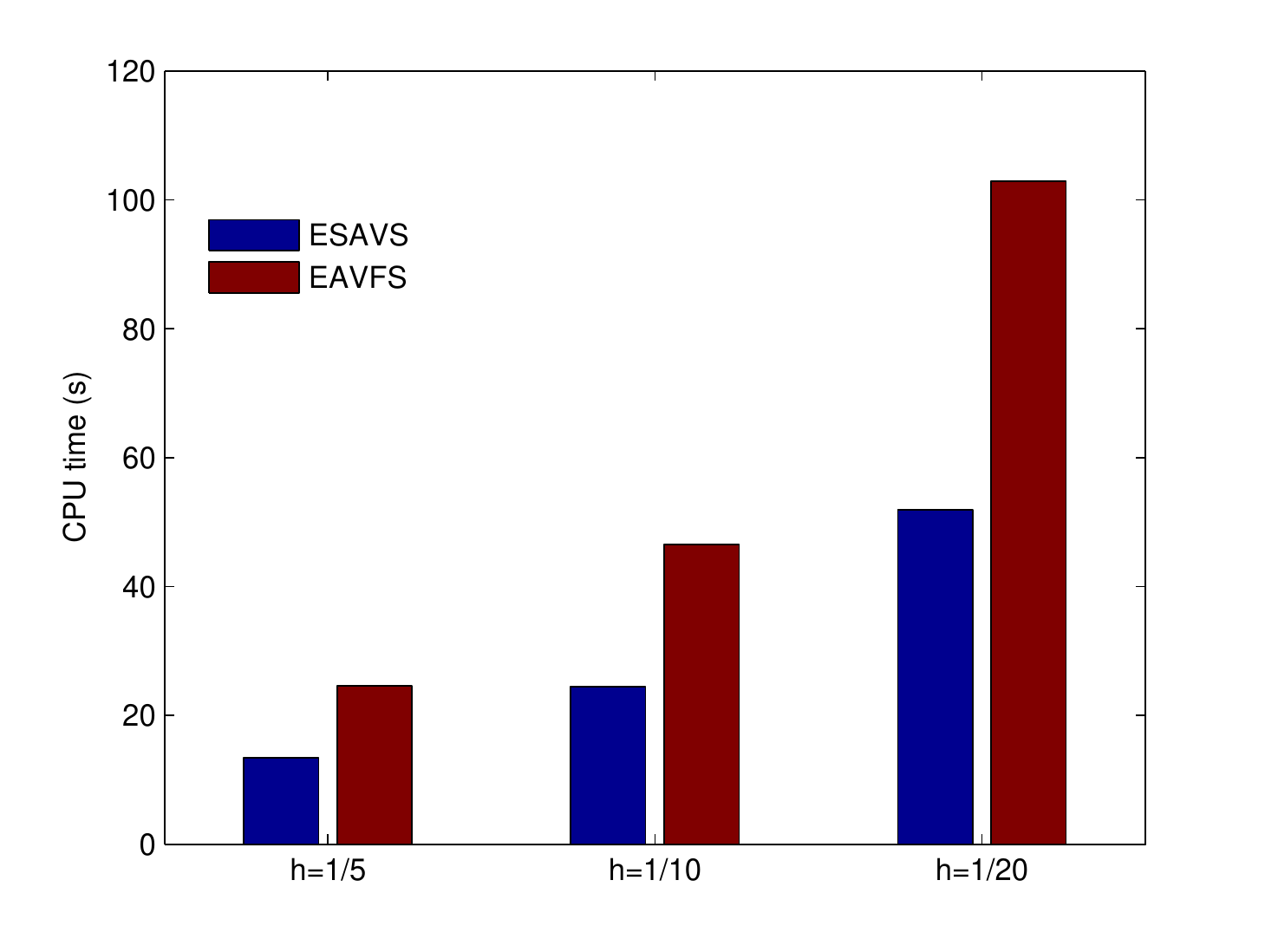}
\caption{CPU time of the two
schemes for the soliton with
different mesh sizes till $t=10$
under $\tau=0.001$.  }\label{Fig_NKGEs:1}
\end{minipage}
\begin{minipage}[t]{65mm}
\includegraphics[width=65mm]{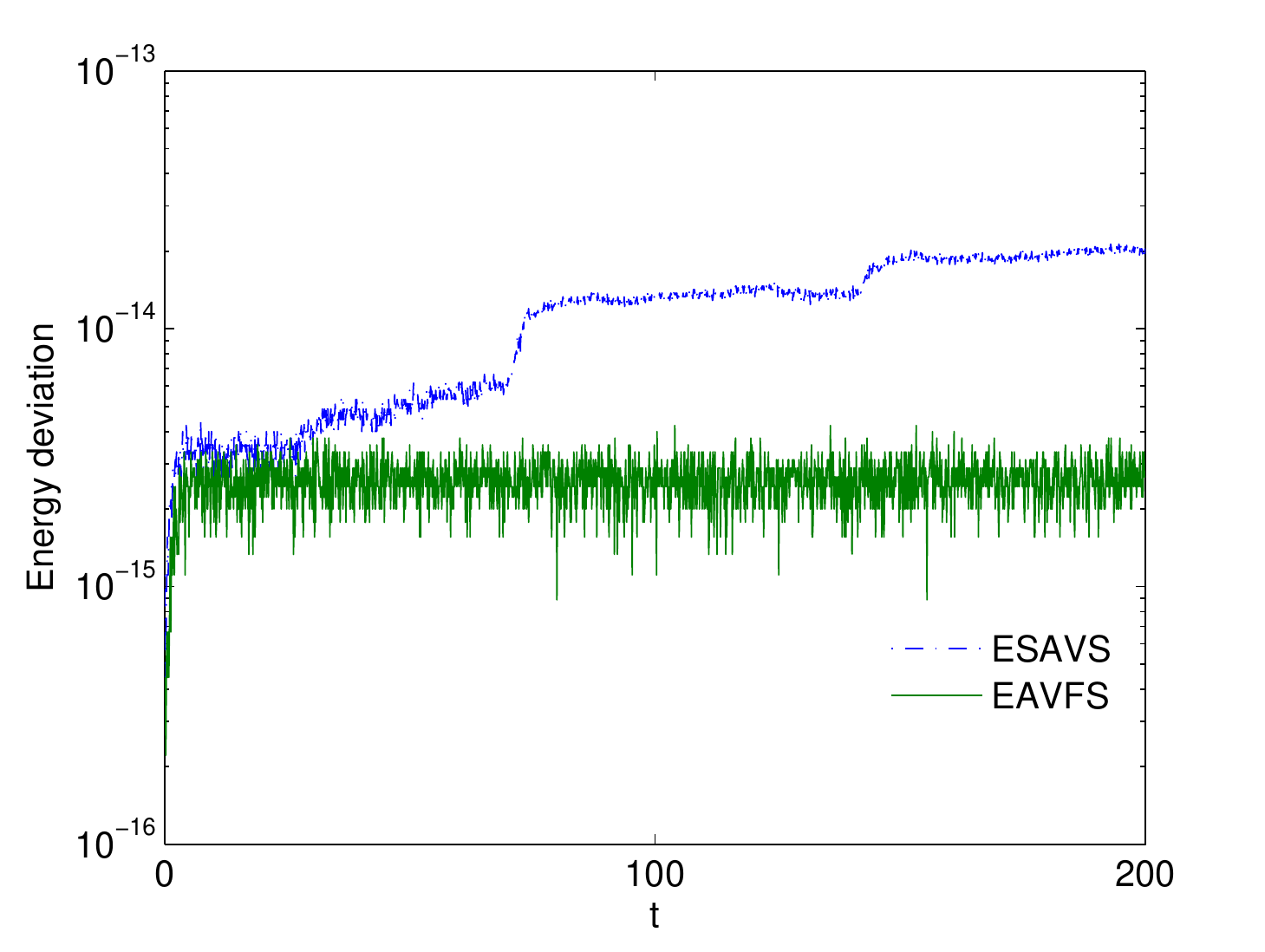}
\caption{The energy deviation with $h=\tau=0.1$ over the time interval $t\in[0,200]$.}\label{Fig_NKGEs:2}
\end{minipage}\
\end{figure}

The error and convergence order of EAVFS and ESAVS at time $t=1$ are given in Tab. \ref{Tab:NKGEs:2}, which can be observed that
all schemes have second order accuracy in time and space and the error provided by ESAVS has the same order of magnitude as the one provided by ESAVS. Besides, we carry out comparison on the computational cost of the two schemes in Fig. \ref{Fig_NKGEs:1} by refining the mesh size gradually, which shows that
the cost of ESAVS is cheaper. Moreover, as the refinement of mesh sizes, the advantage of ESAVS emerges, which implies that our scheme shows the remarkable performance in the efficiency. The long-term energy deviations are plotted in Fig. \ref{Fig_NKGEs:2}. It is clear that  ESAVS and EAVFS can exactly preserve the discrete energies.

\begin{table}[H]
\tabcolsep=9pt
\footnotesize
\renewcommand\arraystretch{1.1}
\centering
\caption{{Numerical error and convergence rate for the two schemes under different grid steps at $t=1$.}}\label{Tab:NKGEs:2}
\begin{tabular*}{\textwidth}[h]{@{\extracolsep{\fill}} c c c c c c}\hline
{Scheme\ \ } &{$(h,\tau)$} &{$L^2$-error} &{order}& {$L^{\infty}$-error}&{order}  \\     
\hline
 \multirow{4}{*}{ESAVS}  &{$(\frac{1}{10},\frac{1}{100})$}& {1.287e-03}&{-} &{1.367e-03} & {-}\\[1ex]
  {}  &{$(\frac{1}{20},\frac{1}{200})$}& {3.217e-04}&{2.00} &{3.413e-04} & {2.00}\\[1ex] 
   {}  &{$(\frac{1}{40},\frac{1}{400})$}& {8.044e-05}&{2.00} &{8.531e-05} &{2.00} \\[1ex]
    {}  &{$(\frac{1}{80},\frac{1}{800})$}& {2.011e-05}&{2.00} &{2.133e-05} &{2.00} \\\hline
 \multirow{4}{*}{EAVFS}  &{$(\frac{1}{10},\frac{1}{100})$}& {1.104e-03}&{-} &{1.050e-03} & {-}\\[1ex]
  {}  &{$(\frac{1}{20},\frac{1}{200})$}& {2.761e-04}&{2.00} &{2.621e-04} & {2.00}\\[1ex] 
   {}  &{$(\frac{1}{40},\frac{1}{400})$}& {6.902e-05}&{2.00} &{6.551e-05} &{2.00} \\[1ex]
    {}  &{$(\frac{1}{80},\frac{1}{800})$}& {1.725e-05}&{2.00} &{1.638e-05} &{2.00} \\\hline
\end{tabular*}
\end{table}

Then, we apply the proposed scheme to solve the following two dimensional nonlinear sine-Gordon equation
\begin{align*}
&\partial_{tt}u-\partial_{xx} u-\partial_{yy} u+\sin(u)=0,\ (x,y)\in\mathbb{R}^2,\ t>0,
\end{align*}
with initial conditions \cite{DPT95}
\begin{align*}
&u(x,y,0)=4\tan^{-1}\left[\exp\left(\frac{4-\sqrt{(x+3)^{2}+(y+7)^{2}}}{0.436}\right)\right],\\
& u_t(x,y,0)=\frac{4.13}{\cosh\left(\frac{4-\sqrt{(x+3)^{2}+(y+7)^{2}}}{0.436}\right)},\ (x,y)\in\mathbb{R}^2.
\end{align*}

\begin{figure}[H]
\centering
\begin{minipage}[t]{70mm}
\includegraphics[width=70mm]{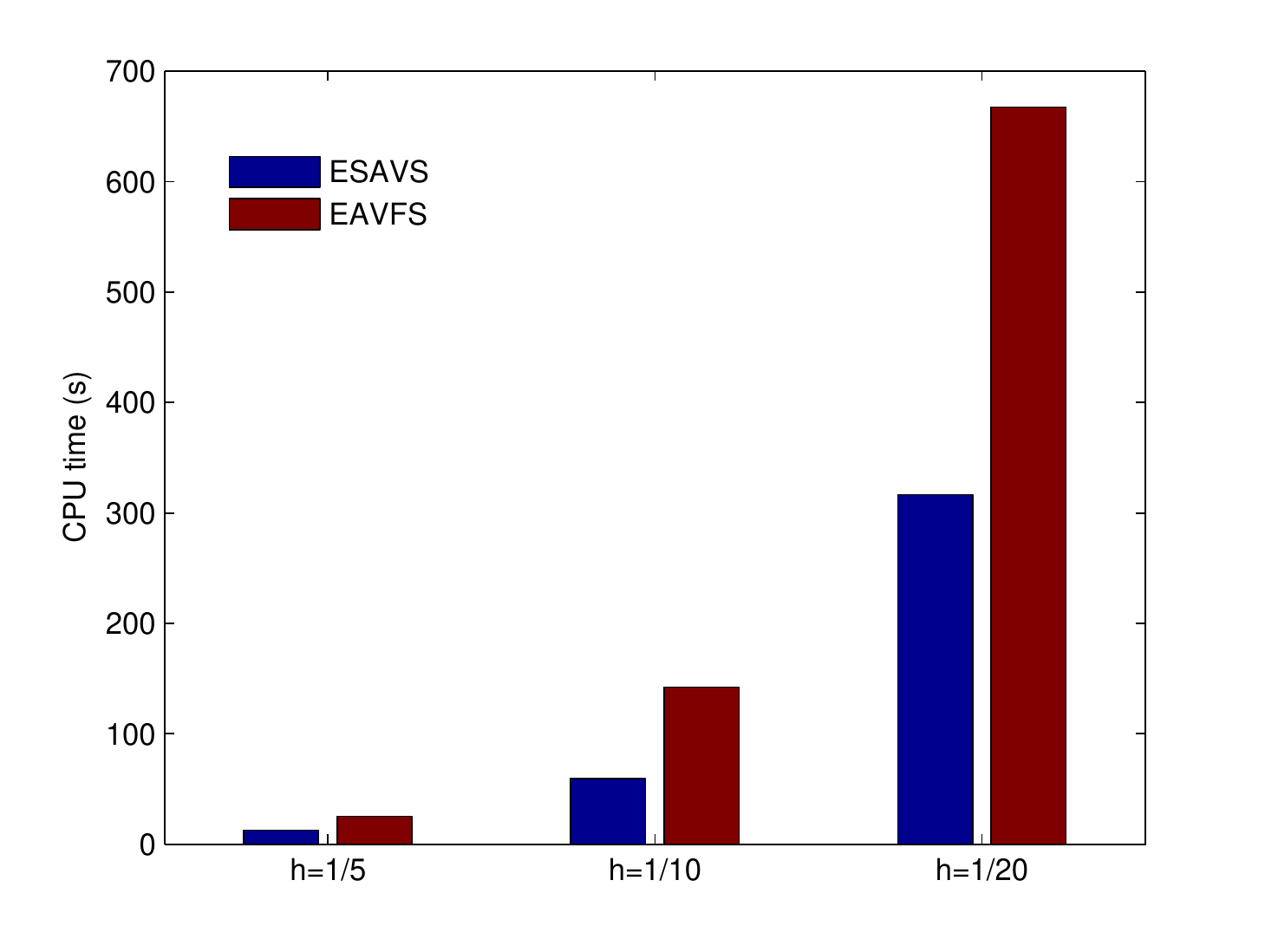}
\end{minipage}
\caption{CPU time of the two
schemes for the soliton with
different mesh sizes till $t=1$
under $\tau=0.01$.}\label{Fig_NKGEs-cpu}
\end{figure}

\begin{figure}[H]
\begin{minipage}[t]{60mm}
\includegraphics[width=60mm]{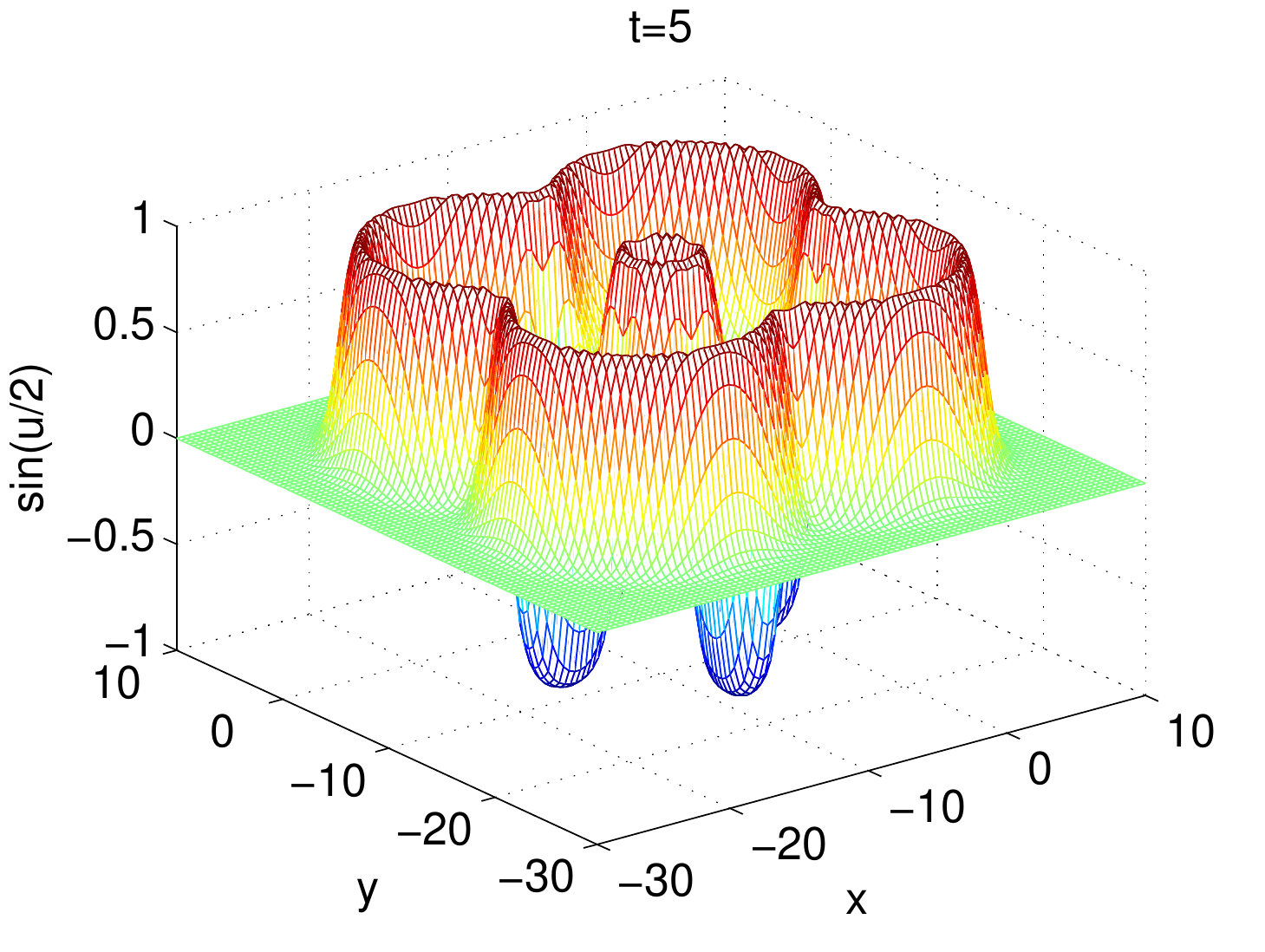}
\end{minipage}
\begin{minipage}[t]{60mm}
\includegraphics[width=60mm]{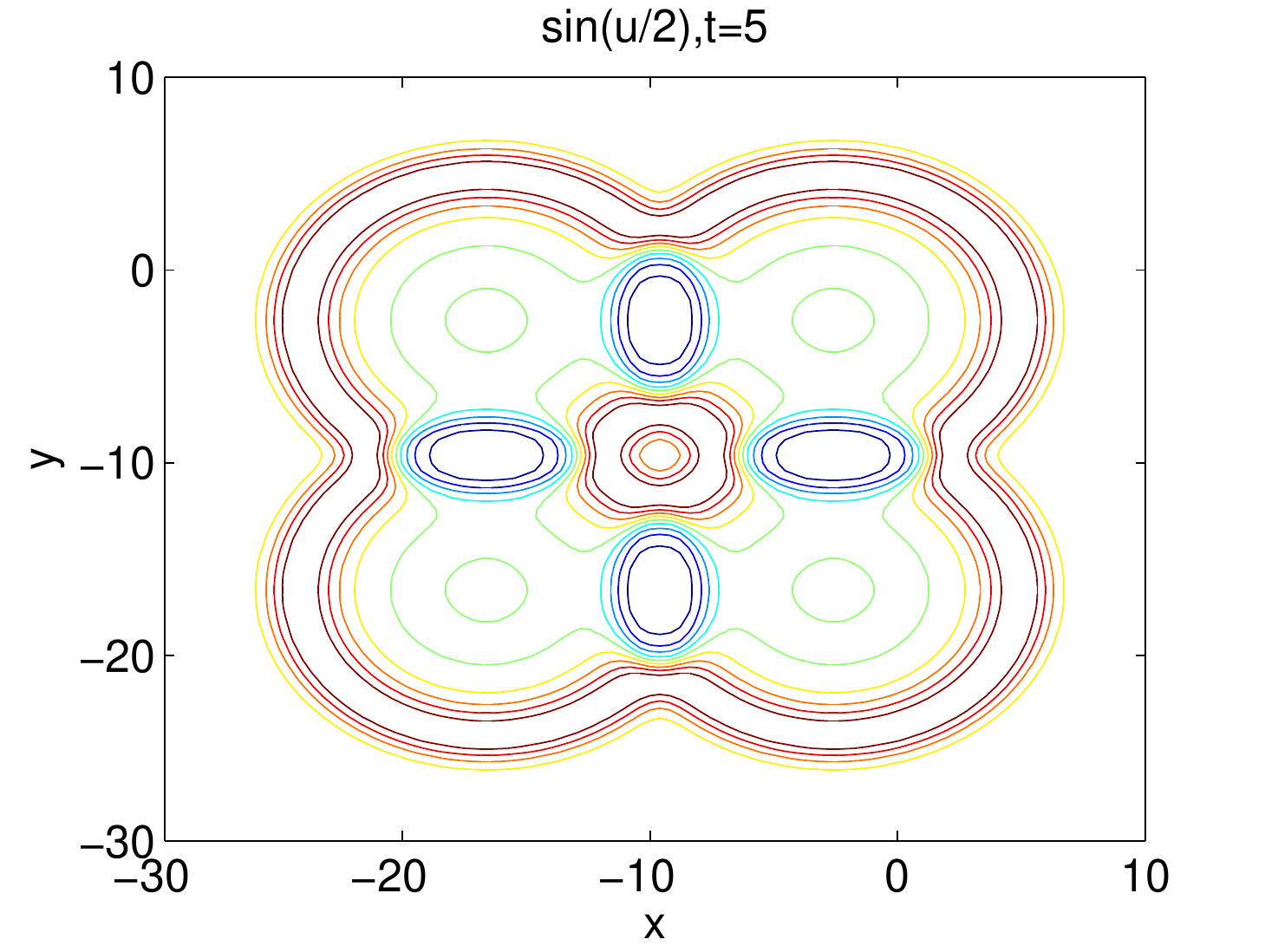}
\end{minipage}
\centering \begin{minipage}[t]{60mm}
\includegraphics[width=60mm]{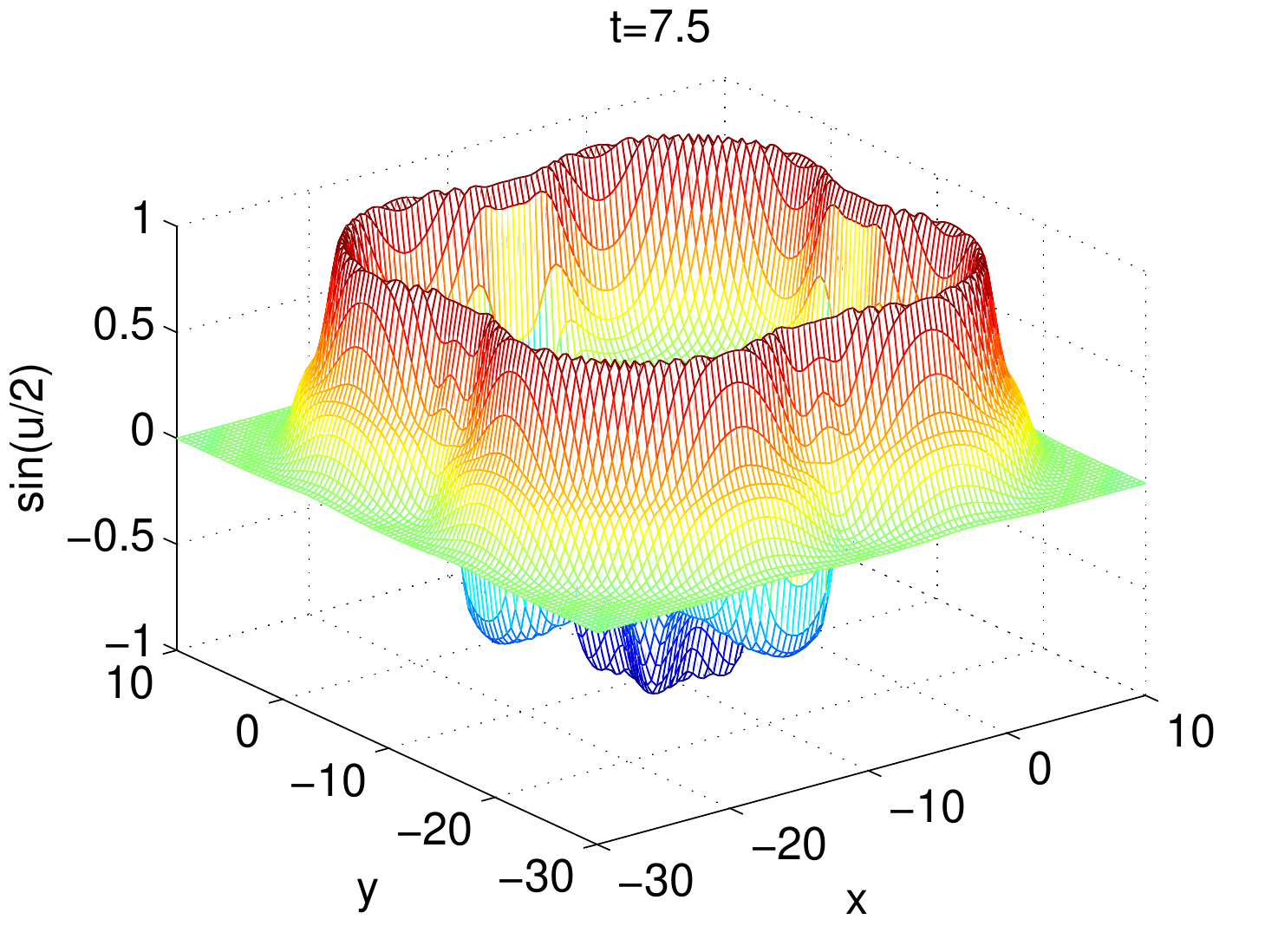}
\end{minipage}
\begin{minipage}[t]{60mm}
\includegraphics[width=60mm]{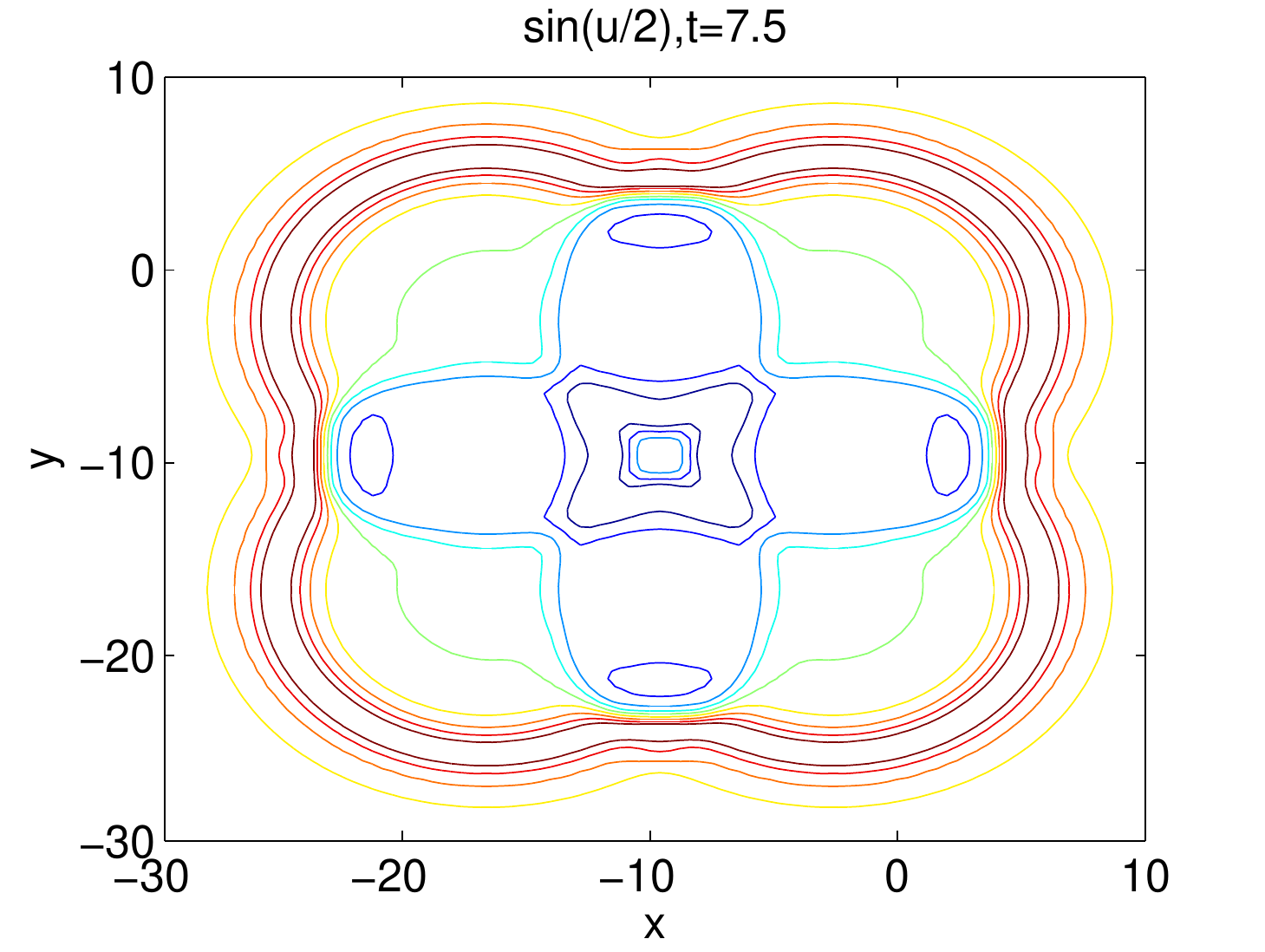}
\end{minipage}
\begin{minipage}[t]{60mm}
\includegraphics[width=60mm]{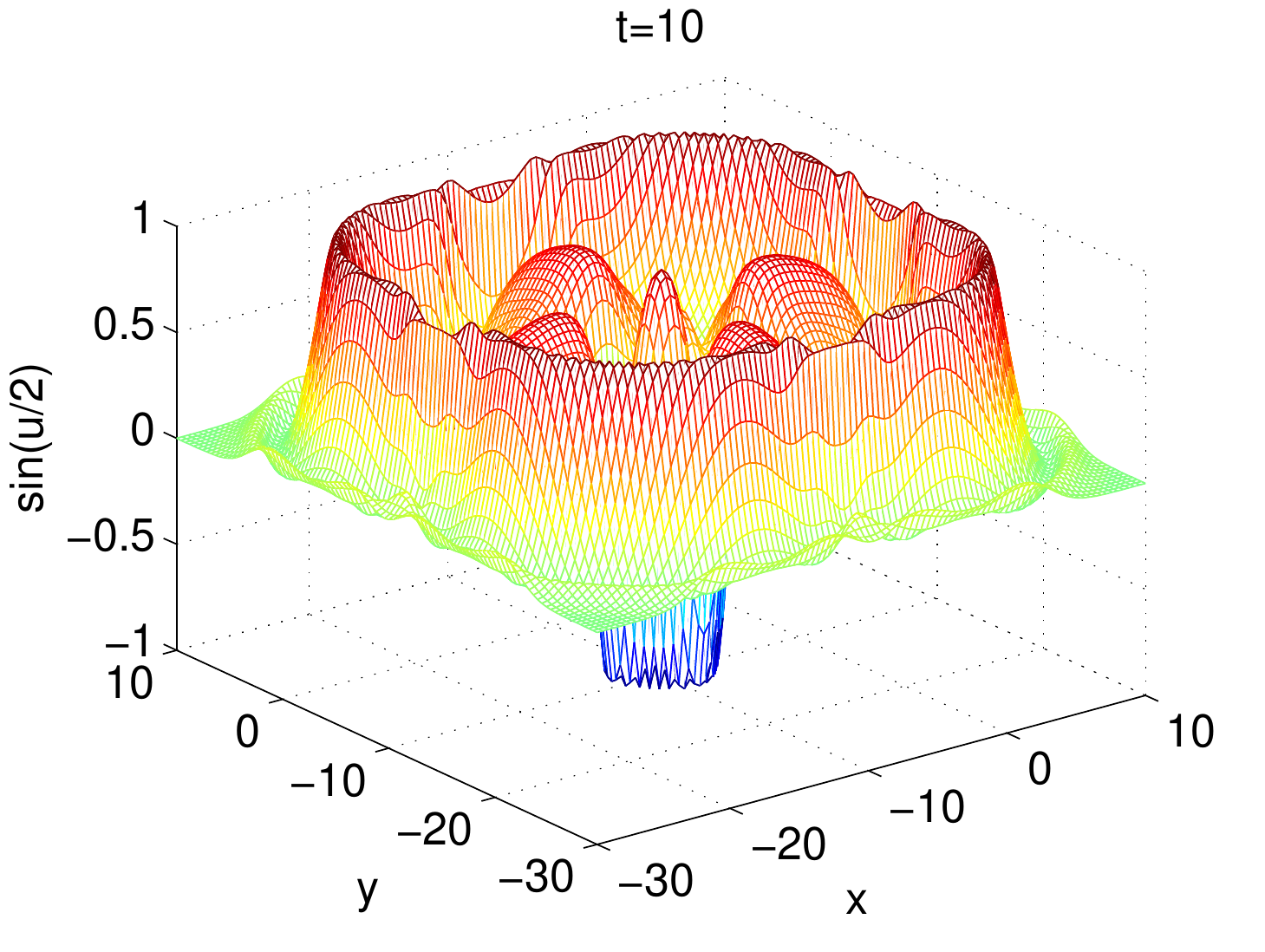}
\end{minipage}
\begin{minipage}[t]{60mm}
\includegraphics[width=60mm]{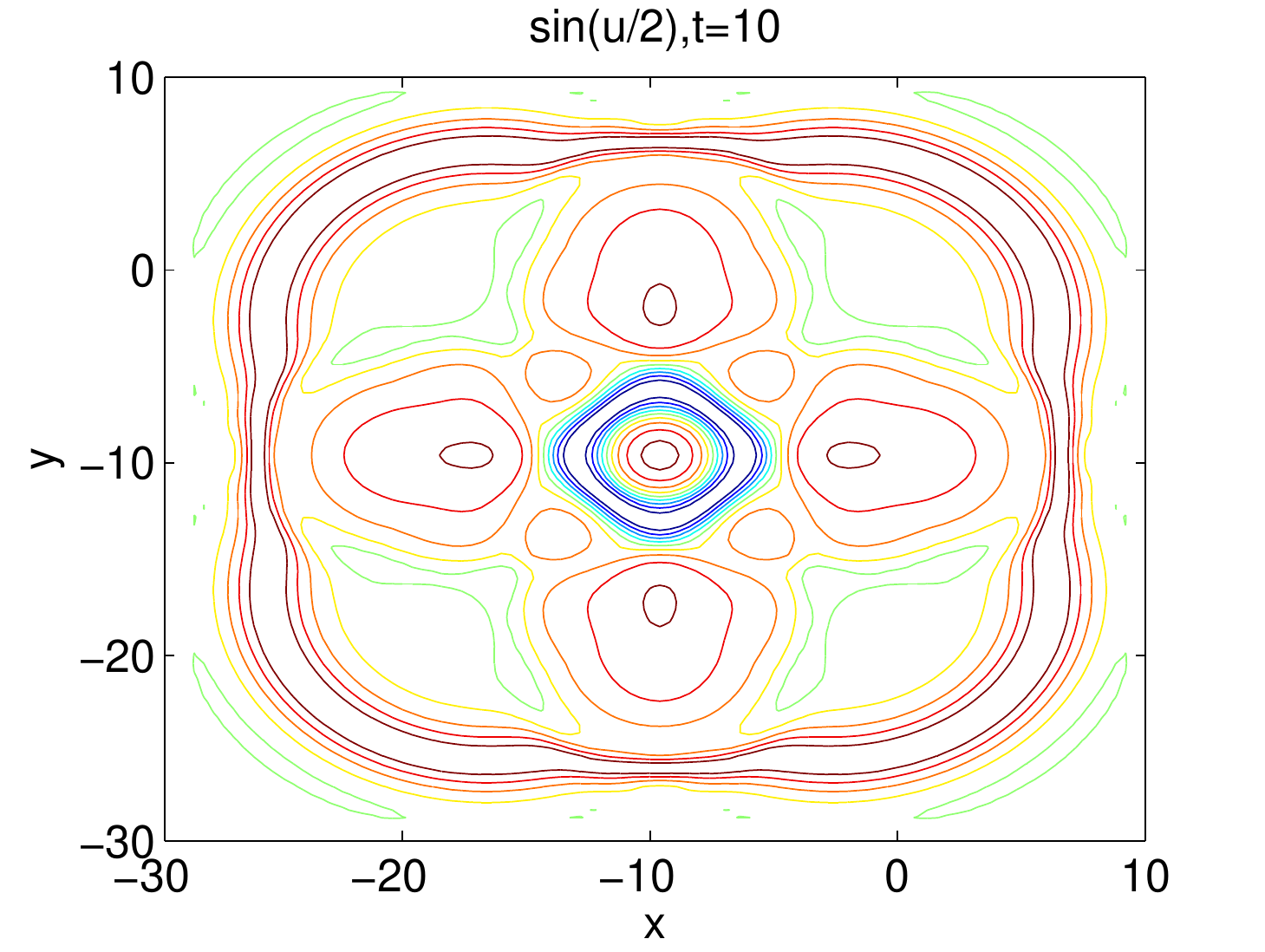}
\end{minipage}
\caption{Collision of four ring solitons (mesh plot (left) and contour plot (right)) in terms of $\sin(u/2)$ at times $t=0,2.5,5,7.5, 10$
with $h=0.2$ and $\tau=0.1$. }\label{KGNEs:3}
\end{figure}

\begin{figure}[H]
\centering\begin{minipage}[t]{70mm}
\includegraphics[width=70mm]{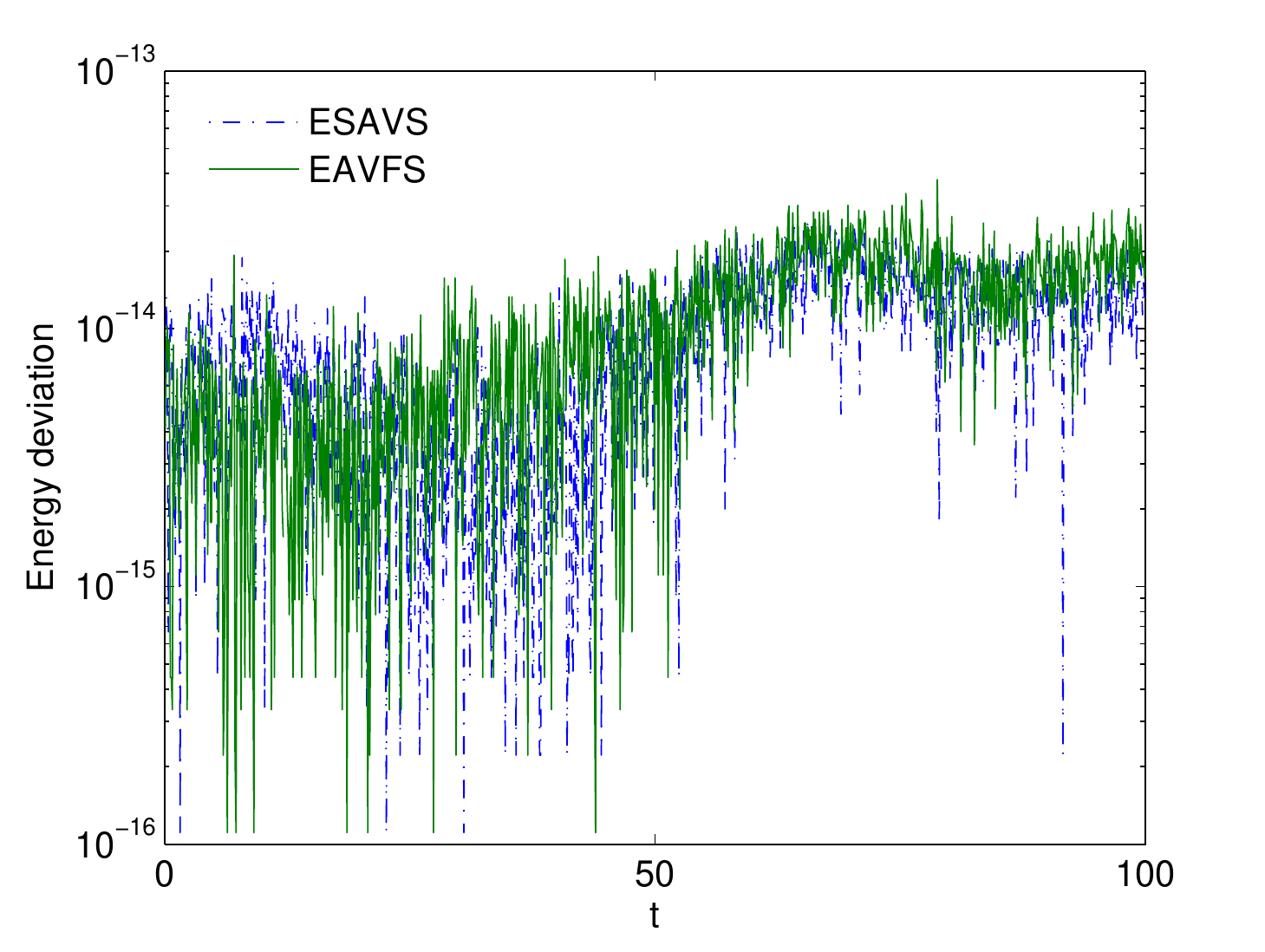}
\end{minipage}
\caption{ The energy deviation over the time interval $t\in[0,100]$
with $h=0.2$ and $\tau=0.1$. }\label{KGNEs:4}
\end{figure}

We take computational domain $\Omega=[-30,10]^2$ with a periodic boundary condition and choose the parameter $C_0=0$. In Fig. \ref{Fig_NKGEs-cpu}, we carry out comparison on the computational cost between two schemes by refining the mesh size gradually. It is clear to see that the cost of EAVFS is more expensive. Moreover, as the refinement of mesh sizes, the advantage of ESAVS emerges, which implies that our scheme is more preferable for large scale simulations than the EAVFS.
Fig. \ref{KGNEs:3} shows the collision precisely among four expanding circular ring solitons which are in good agreement with those given in Refs. \cite{DPT95,SKV10}. {Here, we should note that, following Refs. \cite{DPT95,SKV10}, the solution includes the extension across $x=-10$ and $y=-10$ by symmetry properties of the problem, and the numerical solution in terms of $\sin(u/2)$ instead of $u$ is displayed.} Moverover, we also calculate the energy deviation for the two schemes over the time interval $t\in[0,100]$ and plot it in Fig. \ref{KGNEs:4}. As is clear, our scheme is comparable with the EAVFS.

\begin{figure}[H]
\centering\begin{minipage}[t]{60mm}
\includegraphics[width=60mm]{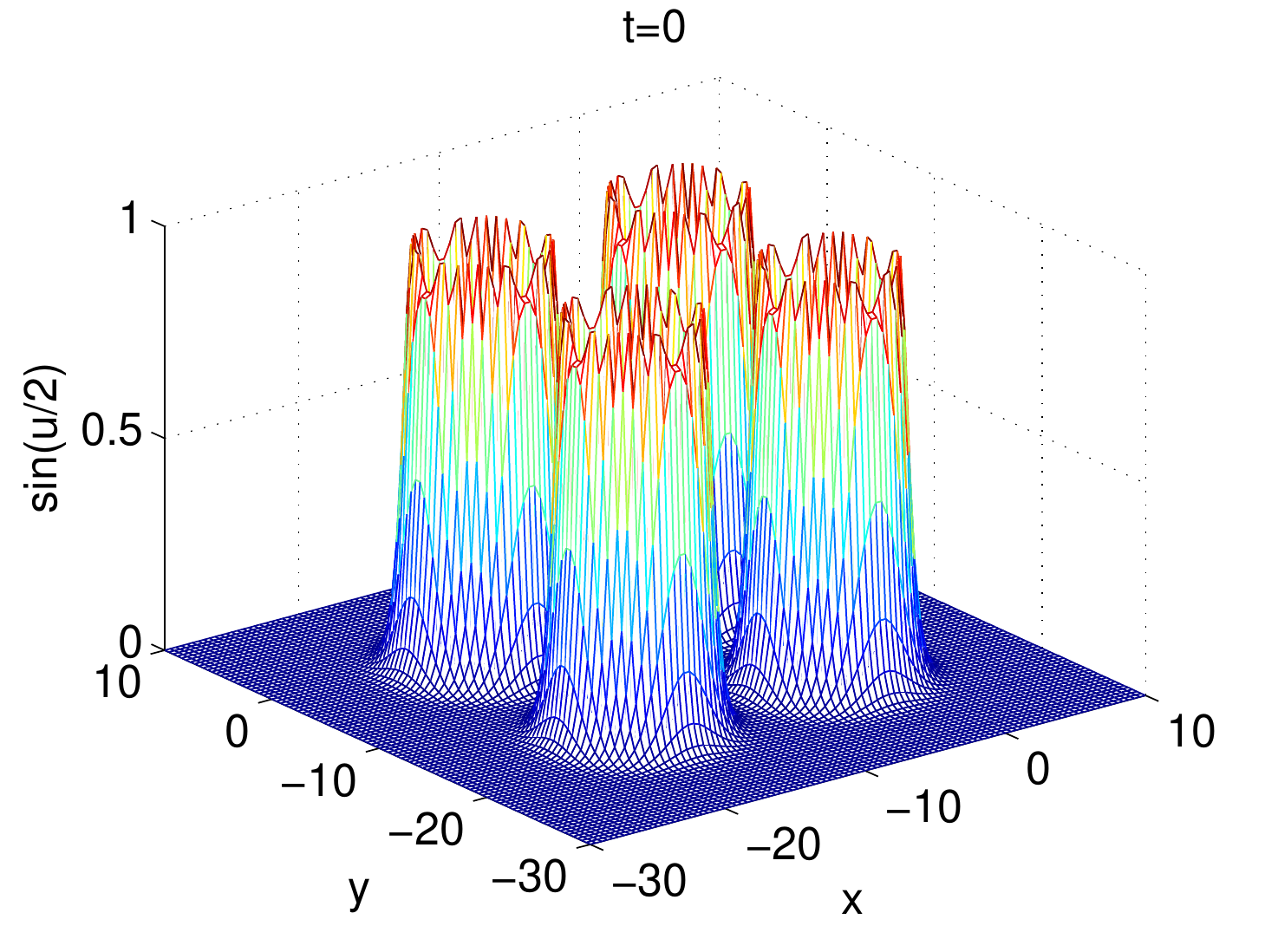}
\end{minipage}
\begin{minipage}[t]{60mm}
\includegraphics[width=60mm]{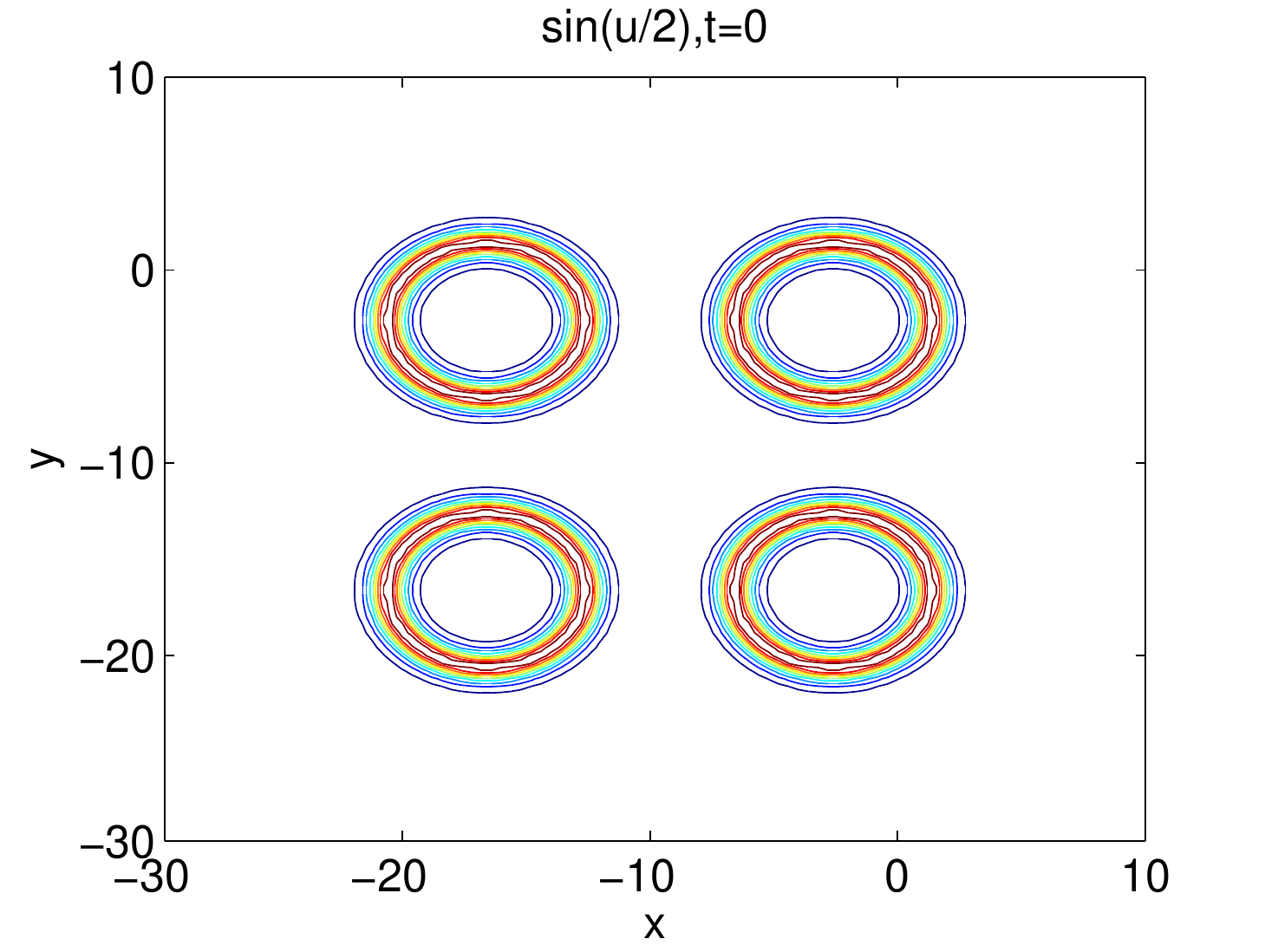}
\end{minipage}
\begin{minipage}[t]{60mm}
\includegraphics[width=60mm]{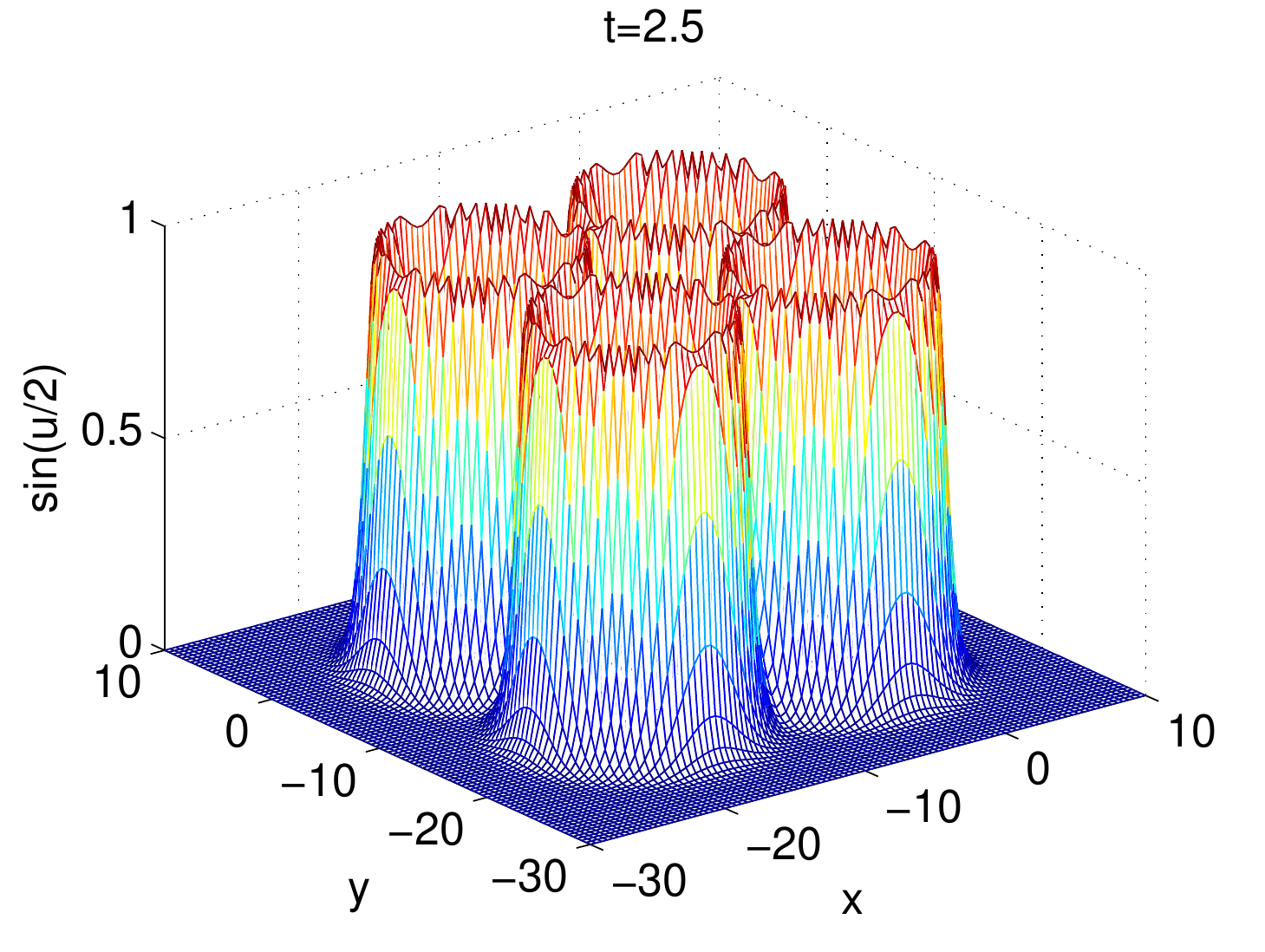}
\end{minipage}
\begin{minipage}[t]{60mm}
\includegraphics[width=60mm]{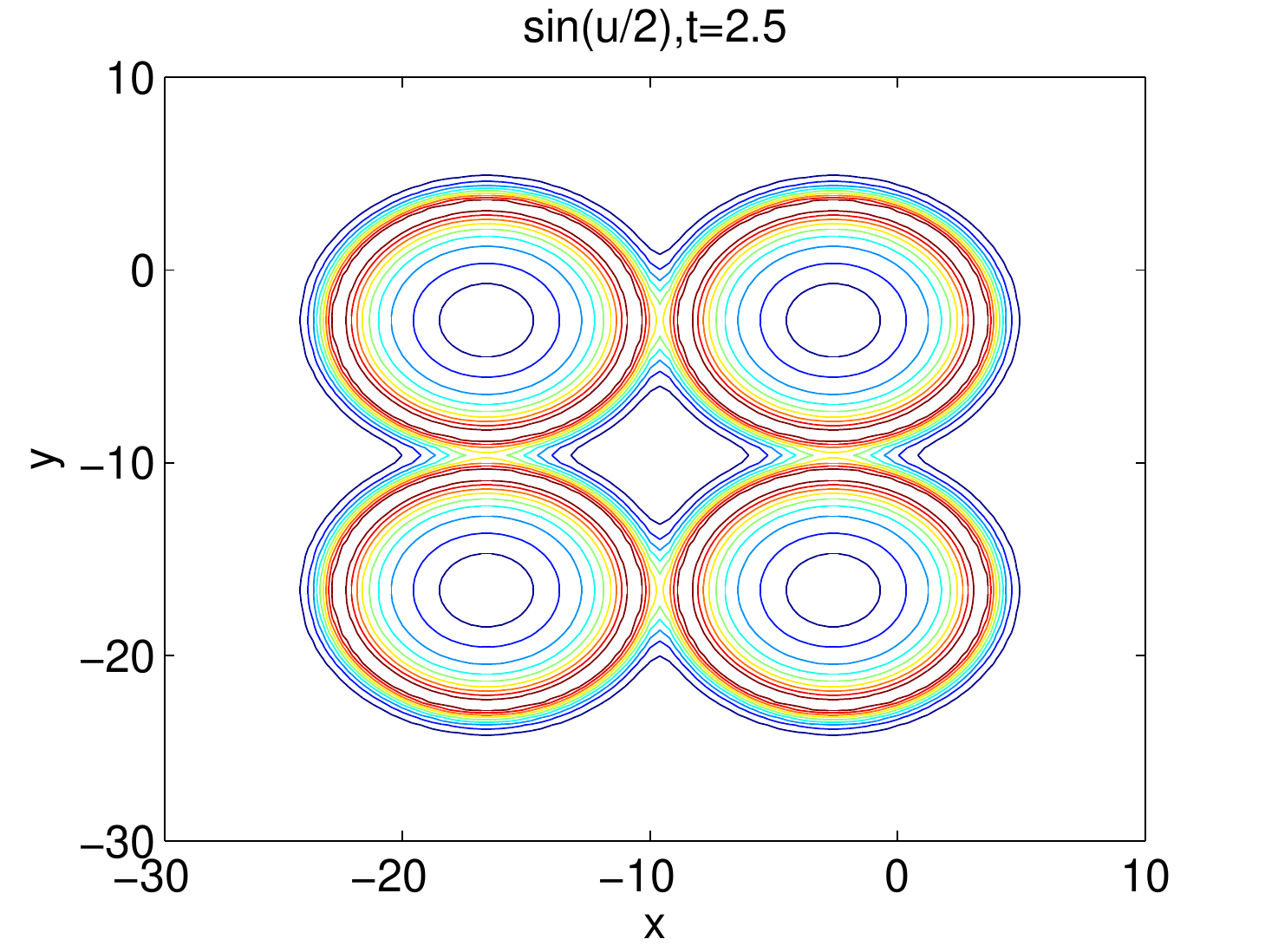}
\end{minipage}
\end{figure}

Finally, we consider the two dimensional nonlinear Klein-Gordon equation, as follows
\begin{align*}
&\partial_{tt}u-\partial_{xx}u-\partial_{yy}u+u^3=0,\ (x,y)\in\mathbb{R}^2,\ t>0,
\end{align*}
with initial conditions 
\begin{align*}
u(x,y,0)=2\text{sech}(\cosh(x^2+y^2)),\ u_t(x,y,0)=0,\ (x,y)\in\mathbb{R}^2.
\end{align*}

\begin{figure}[H]
\centering\begin{minipage}[t]{45mm}
\includegraphics[width=45mm]{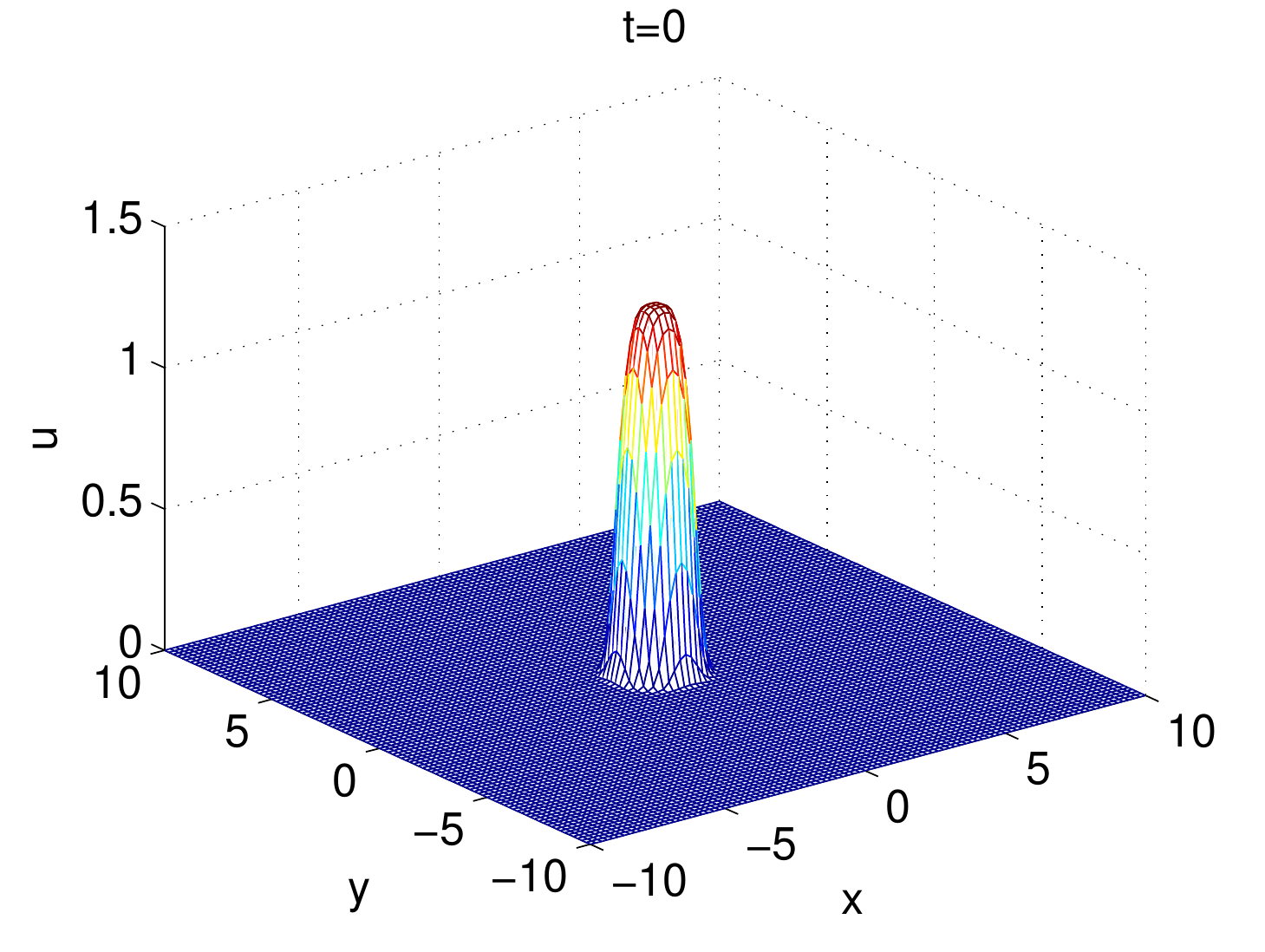}
\end{minipage}
\begin{minipage}[t]{45mm}
\includegraphics[width=45mm]{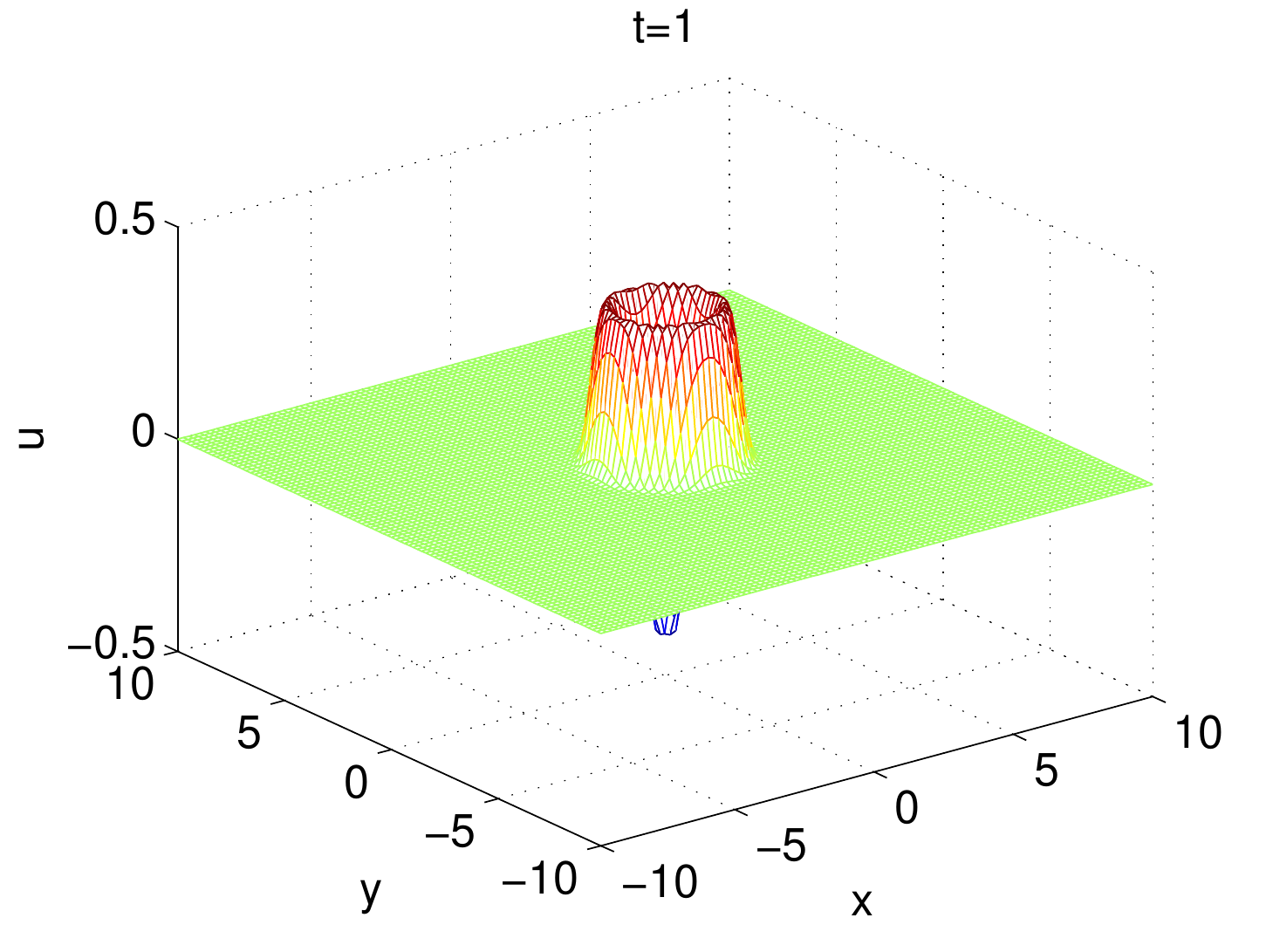}
\end{minipage}
\begin{minipage}[t]{45mm}
\includegraphics[width=45mm]{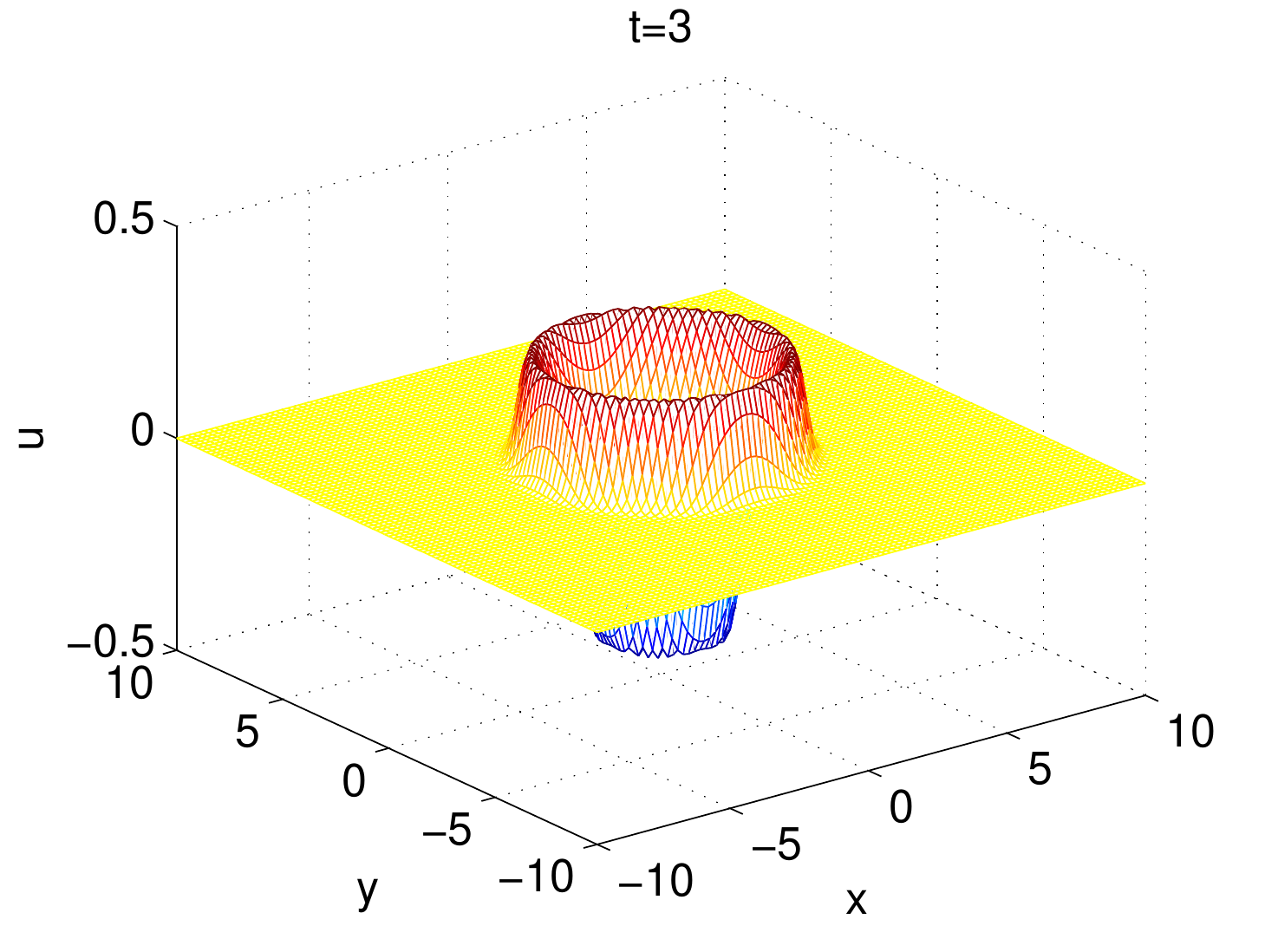}
\end{minipage}
\begin{minipage}[t]{45mm}
\includegraphics[width=45mm]{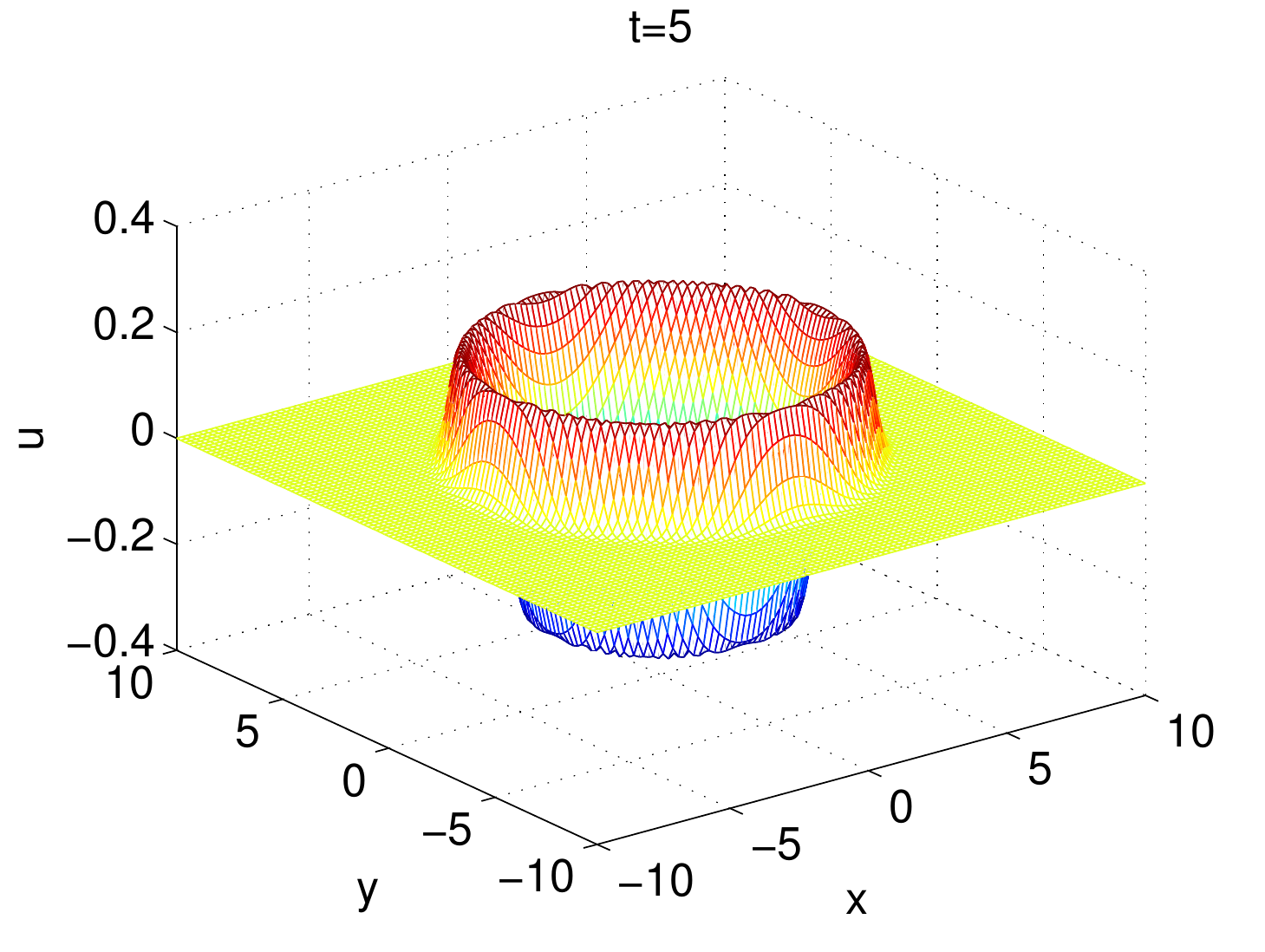}
\end{minipage}
\begin{minipage}[t]{45mm}
\includegraphics[width=45mm]{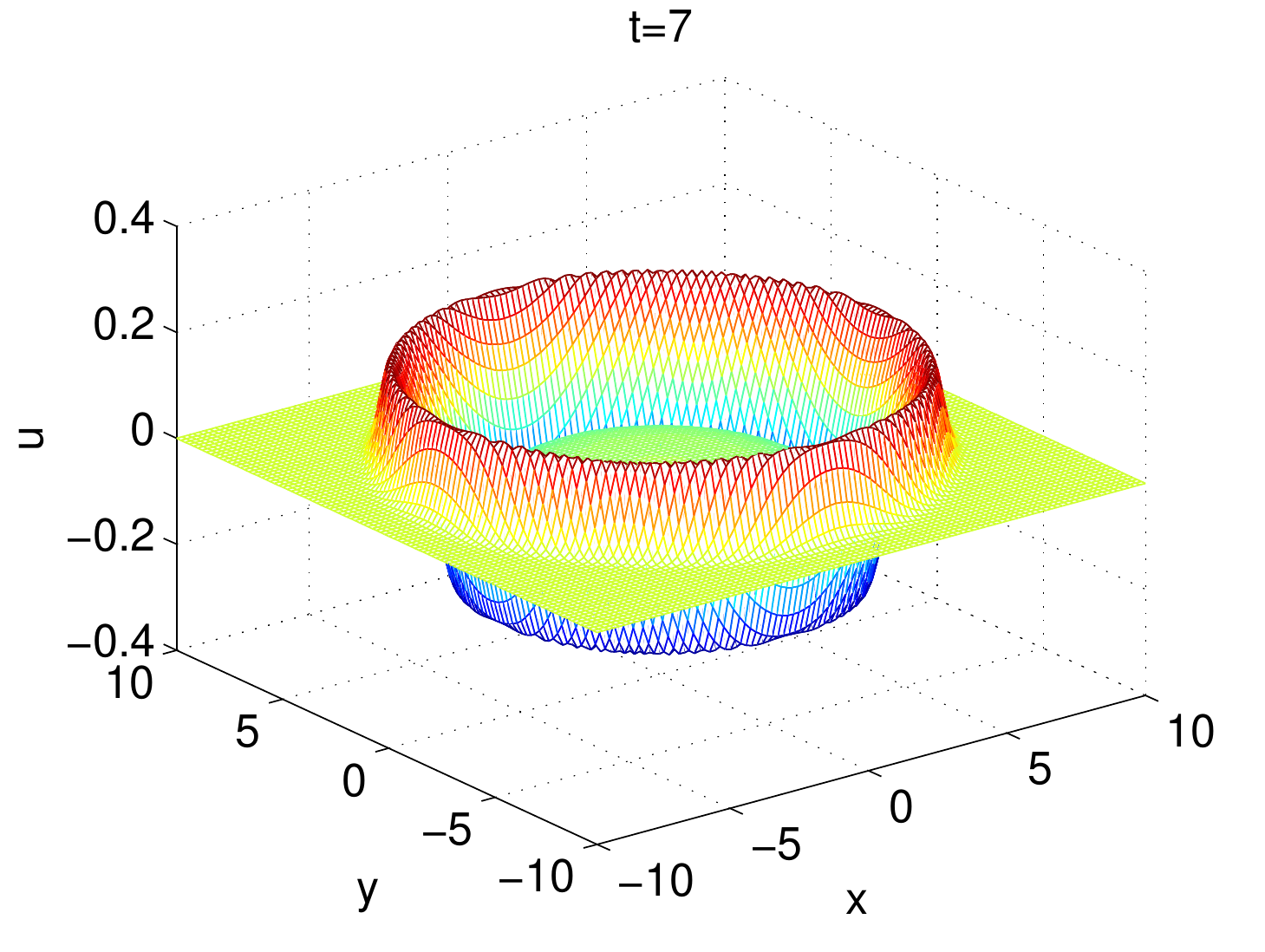}
\end{minipage}
\begin{minipage}[t]{45mm}
\includegraphics[width=45mm]{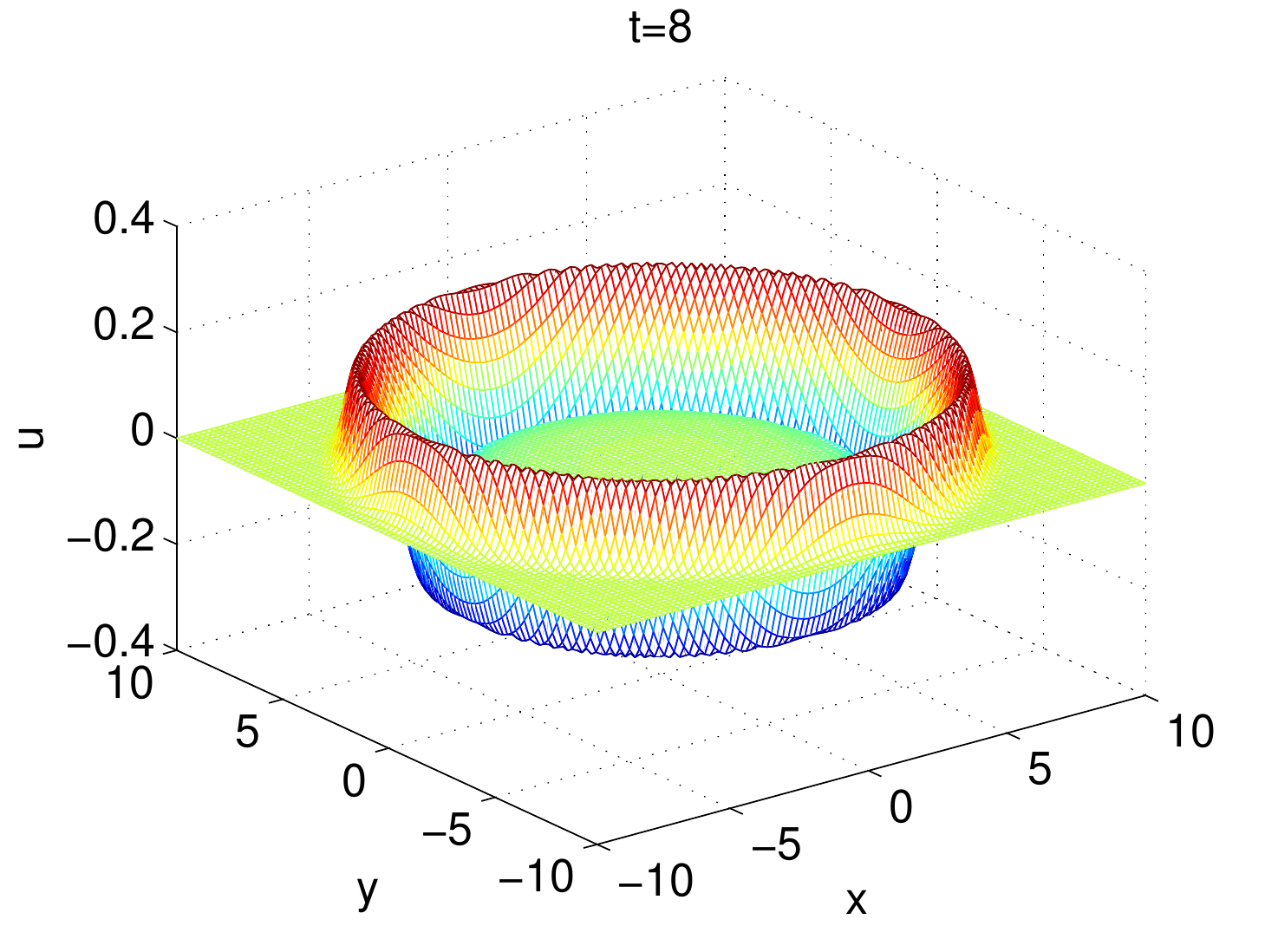}
\end{minipage}
\caption{The Snapshots of numerical solution at times $t=0,1,3,5,7,8$
with $h=\tau=0.1$.}\label{KGNEs:5}
\end{figure}
\begin{figure}[H]
\centering\begin{minipage}[t]{70mm}
\includegraphics[width=70mm]{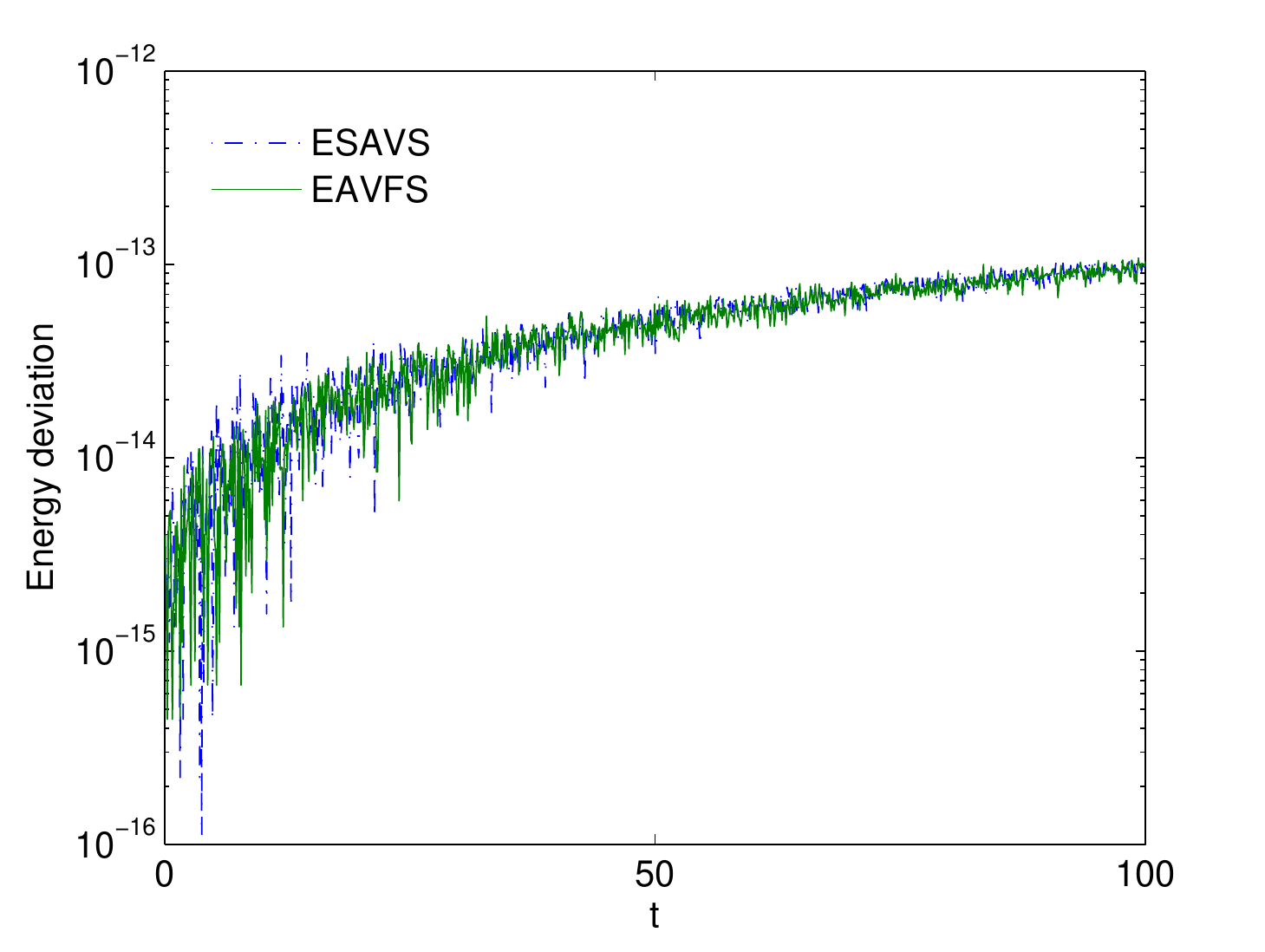}
\end{minipage}
\caption{ The energy deviation over the time interval $t\in[0,100]$
with $h=\tau=0.1$.}\label{KGNEs:6}
\end{figure}

We set the computational domain $\Omega=[-10,10]^2$ with a periodic boundary condition and take the parameter $C_0=0$. Fig. \ref{KGNEs:5} presents the initial condition as well as numerical solutions at different times, which shows the expansion and propagation of the initial soliton to the whole domain until getting the boundary at $t=8$. The long time energy deviation of the two schemes is displayed in Fig. \ref{KGNEs:6}, which behaves similarly as that of Fig. \ref{KGNEs:4}. Here, we omit the comparison of the two schemes for the CPU time. This is because the obtained results behave similarly as that of Fig. \ref{Fig_NKGEs-cpu}.
{\subsection{Nonlinear Schr\"odinger equation}}
In this subsection, we focus on the nonlinear Schr\"odinger equation (NLSE) given as follows
\begin{align}\label{ESAVS-NLS-equation}
\left\lbrace
  \begin{aligned}
&\text{i}\partial_tu({\bm x},t)+\Delta u({\bm x},t)+\beta|u({\bm x},t)|^2u({\bm x},t)=0,\ {\bm x}\in\mathbb{R}^d,\\
&u({\bm x},0)=u_0({\bm x}),\ {\bm x}\in\mathbb{R}^d,
\end{aligned}\right.
  \end{align}
where $\text{i}=\sqrt{-1}$ is the complex unit, $t$ is time variable, ${\bm x}\in\mathbb{R}^d$ is the spatial variable, $u:=u({\bm x},t)$ is the complex-valued wave function, $\Delta$ is the usual Laplace operator, and $\beta$ is a given real constant. The NLSE \eqref{ESAVS-NLS-equation} conserves Hamiltonian energy
\begin{align}\label{ESAVS-NLS-Energy}
E(t):=\int_{\mathbb{R}^d}\Big(-|\nabla u({\bm x},t)|^2+\frac{\beta}{2}|u({\bm x},t)|^4\Big)d{\bm x}=E(0),\ t\ge 0.
\end{align}

For simplicity, we take the one dimensional case i.e. $d=1$ in \eqref{ESAVS-NLS-equation} as example and set the computational domain $\Omega=[a,b]$ with periodic condition boundary. We then let $q:=q(t)=\sqrt{(|u|^4,1)+C_0}$, and rewrite the energy functional \eqref{ESAVS-NLS-Energy} as
\begin{align}\label{ESAVS-NLS-Modified-Energy}
E(t):=\int_{\Omega}\big(-|\nabla u|^2\big)d{x}+\frac{\beta}{2}q^2-\frac{\beta}{2}C_0=E(0),\ t\ge 0.
\end{align}
According to the SAV reformulation, we obtain the following equivalent system
\begin{align}\label{ESAVS-SAV-reformulation}
\left\lbrace
  \begin{aligned}
&\partial_tu=\text{i}\big(\Delta u+\beta \frac{|u|^2u}{\sqrt{(|u|^4,1)+C_0}}q\big),\\
&\partial_t q=\frac{(|u|^2u,\partial_tu)+(\partial_tu,|u|^2u)}{\sqrt{(|u|^4,1)+C_0}},
\end{aligned}\right.
  \end{align}
with the consistent initial condition
\begin{align}\label{ESAVS-SAV-initial}
&u({x},t=0)=u_0({x}),\ q(t=0)=\sqrt{(|u_0({x})|^4,1)+C_0},\ x\in\Omega,
\end{align}
and a periodic condition boundary.

Instead of the finite difference method for discretization of the spatial derivative in \eqref{ESAVS-SAV-reformulation}, we use the standard Fourier pseudo-spectral
method. Actually,  the Laplace operator $\Delta$ is approximated by discrete Fourier transform (DFT) as \begin{align*}
D_2=F^H\Lambda F,
\end{align*}
where
\begin{align*}
\Lambda
=-\Big(\frac{2\pi}{b-a}\Big)^2\text{\rm diag}\Big[0^2,1^2,\cdots,\big(\ \frac{N}{2}\big)^2,\big(-\frac{N}{2}+1\big)^2,\cdots,(-2)^2,(-1)^2\Big].
\end{align*}

 We let $u_{j}^n$ be the numerical approximation of  $u(x_j,t_n)$ for $j=0,1,\cdots,N$ and $n=0,1,2,\cdots$; denote $u^n$ as the solution
vector at $t=t_n$. Then, by an argument similar to the scheme \eqref{NKGEs:eq:3.8}-\eqref{NKGEs:eq:3.9} to system \eqref{ESAVS-SAV-reformulation}, we have

\begin{align}\label{ESAVS-NLS-scheme}
\left\lbrace
  \begin{aligned}
&u^{n+1}=\exp(\text{i}\tau D_2)u^n+\text{i}\beta\tau\int_{0}^1\exp(\text{i}\tau D_2(1-\xi))d\xi \frac{|\hat{u}^{n+\frac{1}{2}}|^2\hat{u}^{n+\frac{1}{2}}q^{n+\frac{1}{2}}}{\sqrt{\langle|\hat{u}^{n+\frac{1}{2}}|^4,{\bm 1}\rangle_{l^2}+C_0}},\\
&q^{n+1}=q^n+\frac{\langle|\hat{u}^{n+\frac{1}{2}}|^2\hat{u}^{n+\frac{1}{2}},u^{n+1}-u^n\rangle_{l^2}+\langle u^{n+1}-u^n,|\hat{u}^{n+\frac{1}{2}}|^2\hat{u}^{n+\frac{1}{2}}\rangle_{l^2}}{\sqrt{\langle|\hat{u}^{n+\frac{1}{2}}|^4,{\bm 1}\rangle_{l^2}+C_0}},
\end{aligned}\right.
  \end{align}
 for $n=1,2,\cdots$. The initial condition \eqref{ESAVS-SAV-initial} is discretized as
\begin{align*}
&u_j^0=u_0(x_j),\ q^0=\sqrt{\langle |u^0|^4,{\bm 1}\rangle_{l^2}+C_0},\ j=0,1,2,\cdots,N.
\end{align*}
\begin{rmk} Note that the proposed scheme \eqref{ESAVS-NLS-scheme} is a three level scheme and we obtain ${u}^1$ and $q^1$ by using ${u}^0$ instead of $\hat{ u}^{\frac{1}{2}}$ for the first step.

  \end{rmk}
\begin{thm}\label{ESAVS-Theorem-5.1} The proposed scheme \eqref{ESAVS-NLS-scheme} preserves the following modified energy
\begin{align*}
{E}_h^{n+1}=E_h^n,\ E_h^n=\langle -D_2u^n,u^n\rangle_{l^2}-\frac{\beta}{2}(q^n)^2-\frac{\beta}{2}C_0.
\end{align*}
\end{thm}
\begin{prf}
The proof is similar to Theorem \ref{ESAV-KNLS-4.1}, thus, for brevity, we omit it.
\end{prf}
Next, we show that the above scheme can be solved efficiently. Eqs. \eqref{ESAVS-NLS-scheme} can be rewritten as
\begin{align}\label{ESAVS-NLS-scheme-eqvi1}
&u^{n+1}=\exp(\text{i}\tau D_2)u^n+\phi\gamma^nq^{n+\frac{1}{2}},\\\label{ESAVS-NLS-scheme-eqvi2}
&q^{n+\frac{1}{2}}=\frac{1}{2}\langle\gamma^n,u^{n+1}\rangle_{l^2}+\frac{1}{2}\langle u^{n+1},\gamma^n\rangle_{l^2}+q^n-\frac{1}{2}\langle\gamma^n,u^{n}\rangle_{l^2}-\frac{1}{2}\langle u^{n},\gamma^n\rangle_{l^2},
\end{align}
where
\begin{align*}
\phi=\text{i}\beta\tau\int_{0}^1\exp(\text{i}\tau D_2(1-\xi))d\xi,\
 \gamma^n=\frac{|\hat{u}^{n+\frac{1}{2}}|^2\hat{u}^{n+\frac{1}{2}}}{\sqrt{\langle|\hat{u}^{n+\frac{1}{2}}|^4,{\bm 1}\rangle_{l^2}+C_0}}.
\end{align*}
Then, by eliminating $q^{n+\frac{1}{2}}$ in \eqref{ESAVS-NLS-scheme-eqvi1}, we have
\begin{align}\label{ESAV-NLS-fast-solver}
u^{n+1}=\frac{1}{2}\phi\gamma^n\langle\gamma^n,u^{n+1}\rangle_{l^2}+\frac{1}{2}\phi\gamma^n\langle u^{n+1},\gamma^n\rangle_{l^2}+b^n,
\end{align}
where
\begin{align*}
b^n=\exp(\text{i}\tau D_2)u^n+\phi\gamma^nq^n-\frac{1}{2}\phi\gamma^n\langle\gamma^n,u^{n}\rangle_{l^2}-\frac{1}{2}\phi\gamma^n\langle u^{n},\gamma^n\rangle_{l^2}.
\end{align*}
We take the inner product of \eqref{ESAV-NLS-fast-solver} with $\gamma^n$ and have, respectively,
\begin{align}\label{ESAV-NLS-fast-solver-1}
&\Big(1-\frac{1}{2}\langle \gamma^n,\phi\gamma^n\rangle_{l^2}\Big)\langle \gamma^n,u^{n+1}\rangle_{l^2}-\frac{1}{2}\langle \gamma^n,\phi\gamma^n\rangle_{l^2}\langle u^{n+1},\gamma^n\rangle_{l^2}=\langle \gamma^n,b^n\rangle_{l^2},\\\label{ESAV-NLS-fast-solver-2}
&-\frac{1}{2}\langle\phi\gamma^n,\gamma^n\rangle_{l^2}\langle\gamma^n,u^{n+1}\rangle_{l^2}+\Big(1-\frac{1}{2}\langle \phi\gamma^n,\gamma^n\rangle_{l^2}\Big)\langle u^{n+1},\gamma^n\rangle_{l^2}=\langle b^n, \gamma^n\rangle_{l^2}.
\end{align}
Eqs. \eqref{ESAV-NLS-fast-solver-1} and \eqref{ESAV-NLS-fast-solver-2} form a $2\times 2$ linear system for the unknowns $\Big( \langle \gamma^n,u^{n+1}\rangle_{l^2},\langle u^{n+1},\gamma^n\rangle_{l^2}\Big)^T$.

Solving $\Big( \langle \gamma^n,u^{n+1}\rangle_{l^2}, \langle u^{n+1},\gamma^n\rangle_{l^2}\Big)^T$ from the $2\times 2$ linear system \eqref{ESAV-NLS-fast-solver-1} and \eqref{ESAV-NLS-fast-solver-2}, and $u^{n+1}$ is then updated from \eqref{ESAV-NLS-fast-solver}. Subsequently, $q^{n+\frac{1}{2}}$ is obtained by \eqref{ESAVS-NLS-scheme-eqvi2}. Finally, we have $q^{n+1}=2q^{n+\frac{1}{2}}-q^n$.

\begin{rmk} In addition, by an argument similar to Remark \ref{ESAV-rk-4.2}, we can deduce that
\begin{align*}
\exp(\text{\rm i}\tau D_2)=F^H \exp(\text{\rm i}\tau\Lambda) F,\ \int_{0}^1\exp(\text{\rm i}\tau D_2(1-\xi))d\xi=F^H\Sigma F,
\end{align*}
where
\begin{align*}
\Sigma=\text{\rm diag}\Big[1,\frac{\exp(\text{\rm i}\tau\lambda_1)-1}{\text{\rm i}\tau\lambda_1},\cdots,\frac{\exp(\text{\rm i}\tau\lambda_{N-1})-1}{\text{\rm i}\tau\lambda_{N-1}}\Big].
\end{align*}
\end{rmk}

We repeat the time step refinement test first and choose the parameter $C_0=0$ and $\beta=2$.
The one dimensional Schr\"odinger equation \eqref{ESAVS-NLS-equation} admits the analytical solution
\begin{align*}
u(x,t)=\text{sech}(x-4t)\exp(2\text{i}x-3\text{i}t),\ x\in \mathbb{R}.
\end{align*}
We choose the analytical solution at $t=0$ as initial condition and set the computational domain $\Omega=[-40,40]$ with a periodic boundary.
To test the temporal discretization errors of the two numerical schemes, we fix the Fourier node $4096$ such that the spatial discretization errors are negligible.

\begin{figure}[H]
\centering\begin{minipage}[t]{60mm}
\includegraphics[width=60mm]{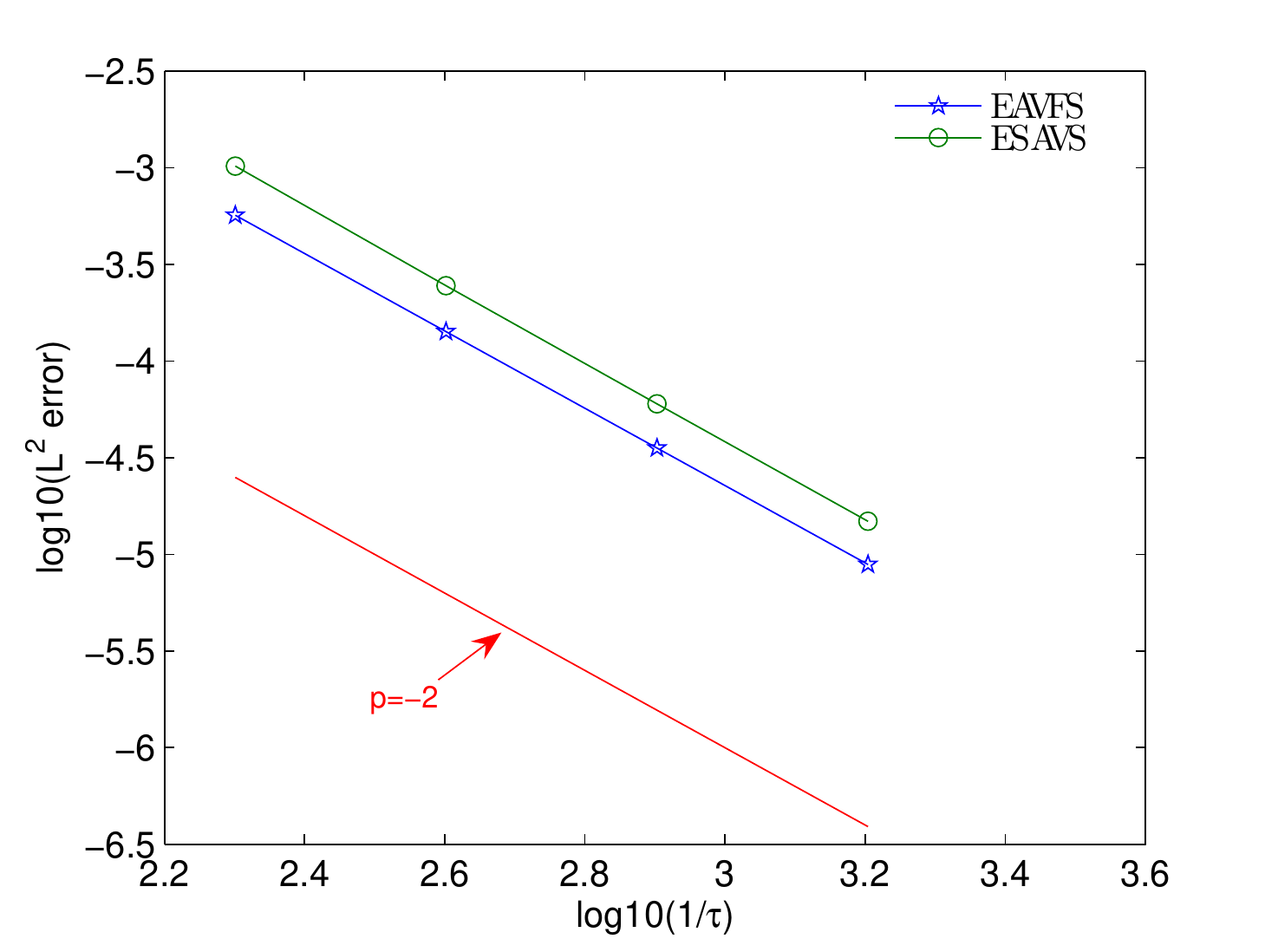}
\end{minipage}
\begin{minipage}[t]{60mm}
\includegraphics[width=60mm]{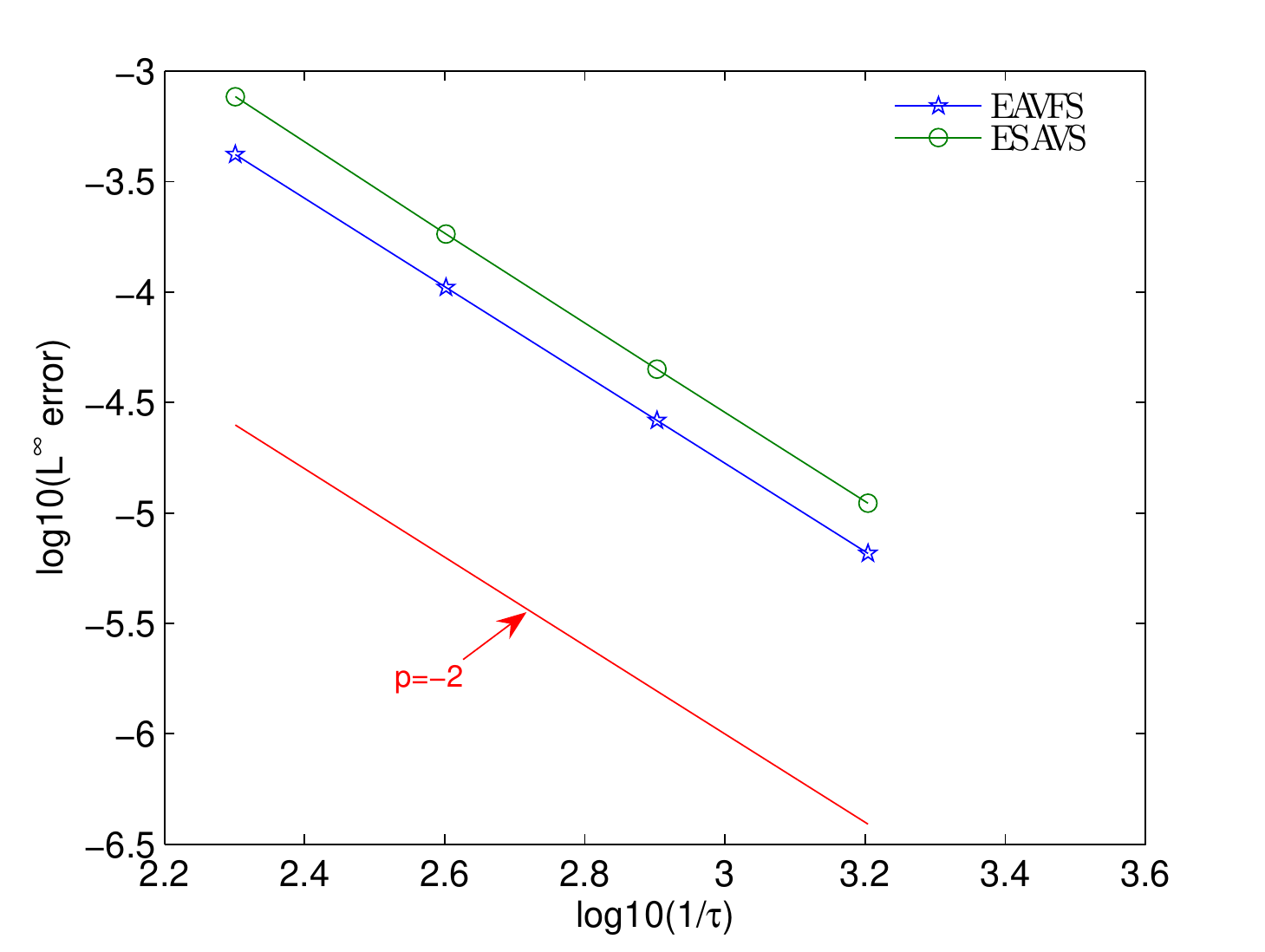}
\end{minipage}
\caption{ Time step refinement tests using the two numerical schemes for the one dimensional Schr\"odinger equation.}\label{ESAV-scheme:fig:1}
\end{figure}

\begin{figure}[H]
\centering\begin{minipage}[t]{60mm}
\includegraphics[width=60mm]{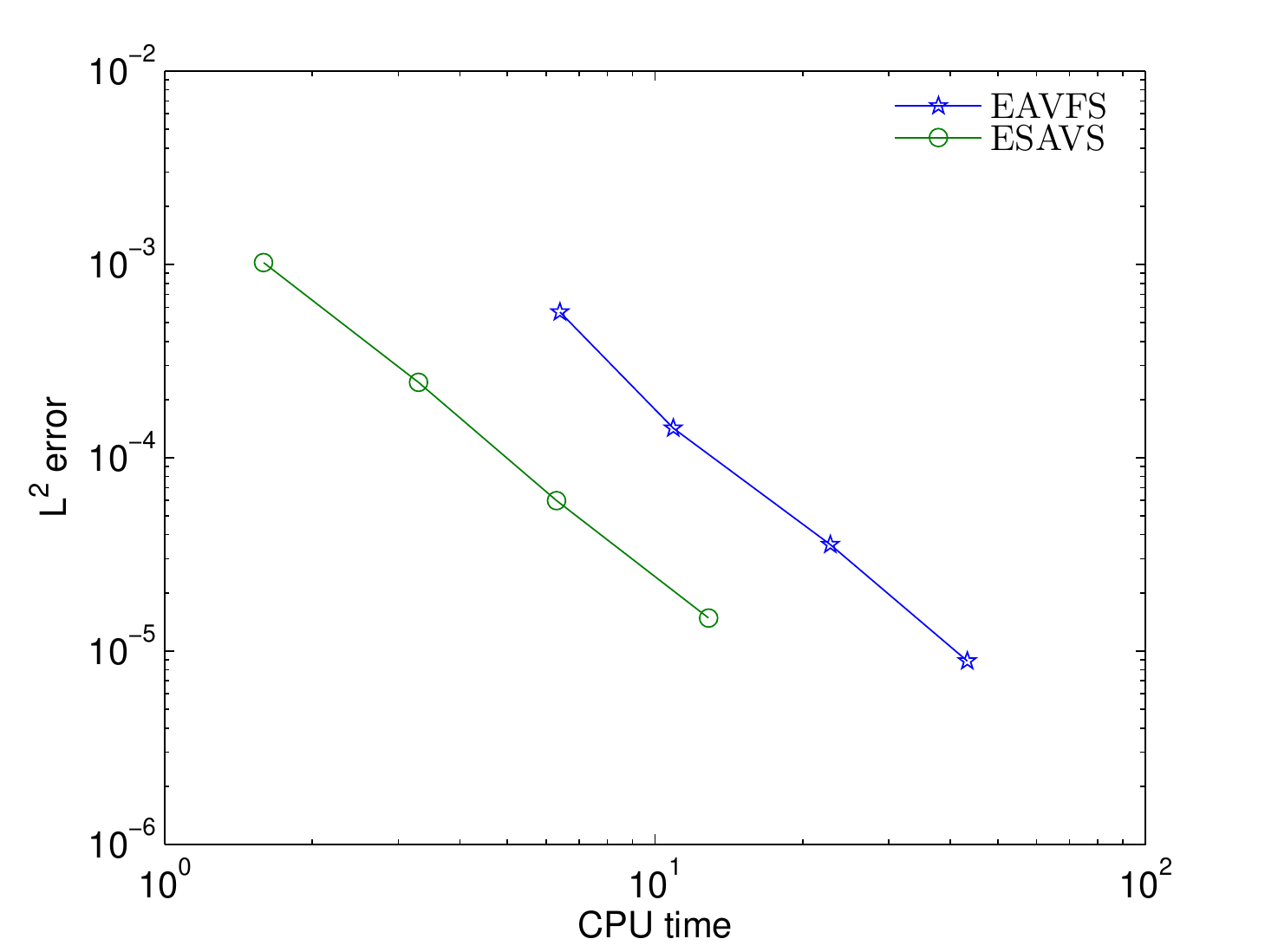}
\end{minipage}
\begin{minipage}[t]{60mm}
\includegraphics[width=60mm]{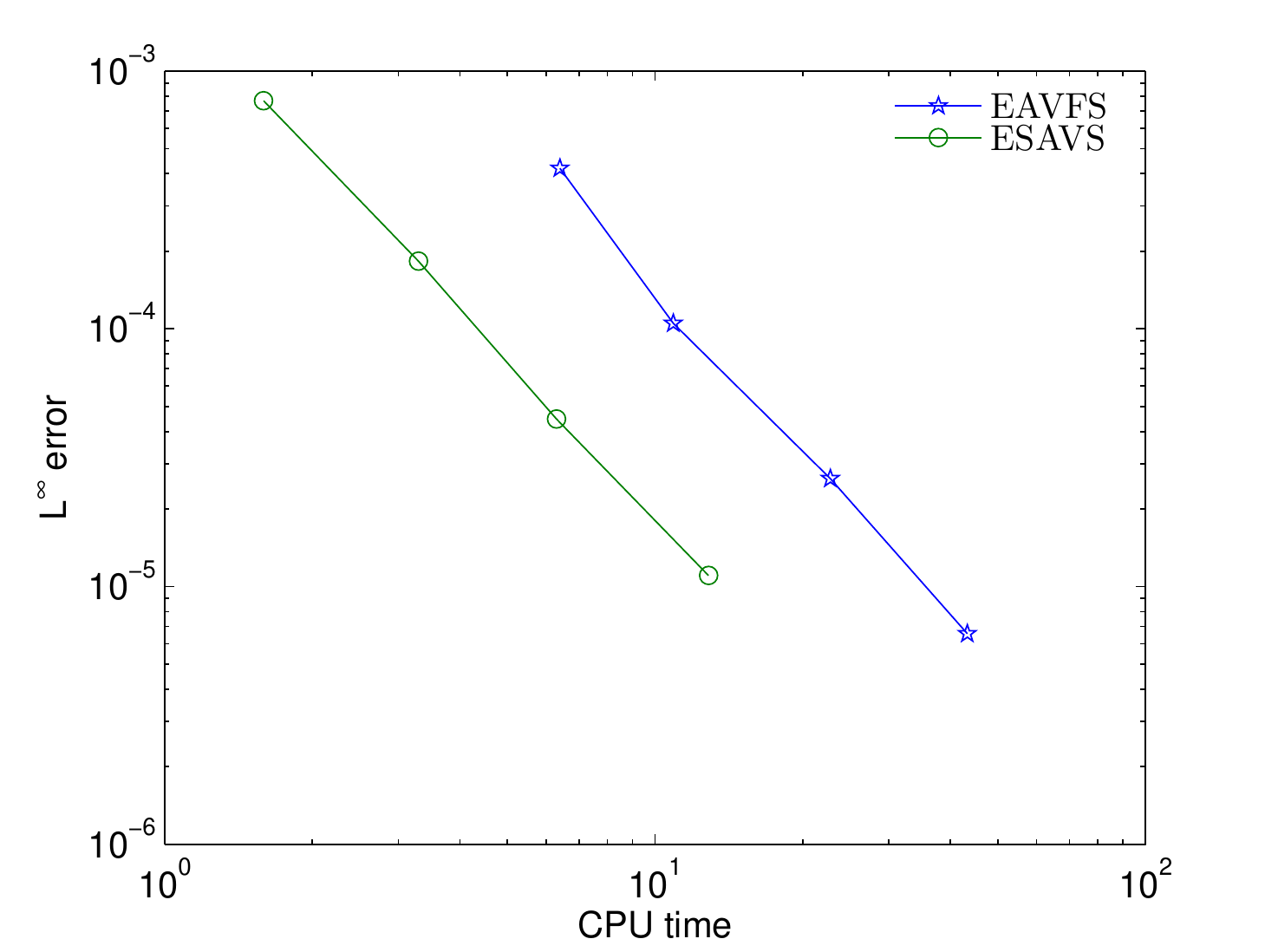}
\end{minipage}
\caption{ The numerical error versus the CPU time using the two numerical schemes for the one dimensional Schr\"odinger equation.}\label{ESAV-scheme:fig:2}
\end{figure}

The $L^2$ errors and $L^{\infty}$ errors in numerical solution of $u$ at $t=1$ are calculated using two numerical schemes with various
time steps, and the results are displayed in Fig. \ref{ESAV-scheme:fig:1}. In Fig. \ref{ESAV-scheme:fig:2},
we show the global $L^2$ errors and $L^{\infty}$ errors of $u$ versus the CPU time using the two different schemes at $t=1$.
From Figs. \ref{ESAV-scheme:fig:1} and \ref{ESAV-scheme:fig:2}, we can draw the following observations:
(i) all schemes have second order accuracy in time; (ii) the error provided by the EAVFS is smallest,
 and the one provided by the proposed scheme has the same order of magnitude as the one of the ESAV scheme;
(iii) for a given global error, the cost of the EAVFS is more expensive than the proposed scheme.

\begin{figure}[H]
\centering
\begin{minipage}[t]{70mm}
\includegraphics[width=70mm]{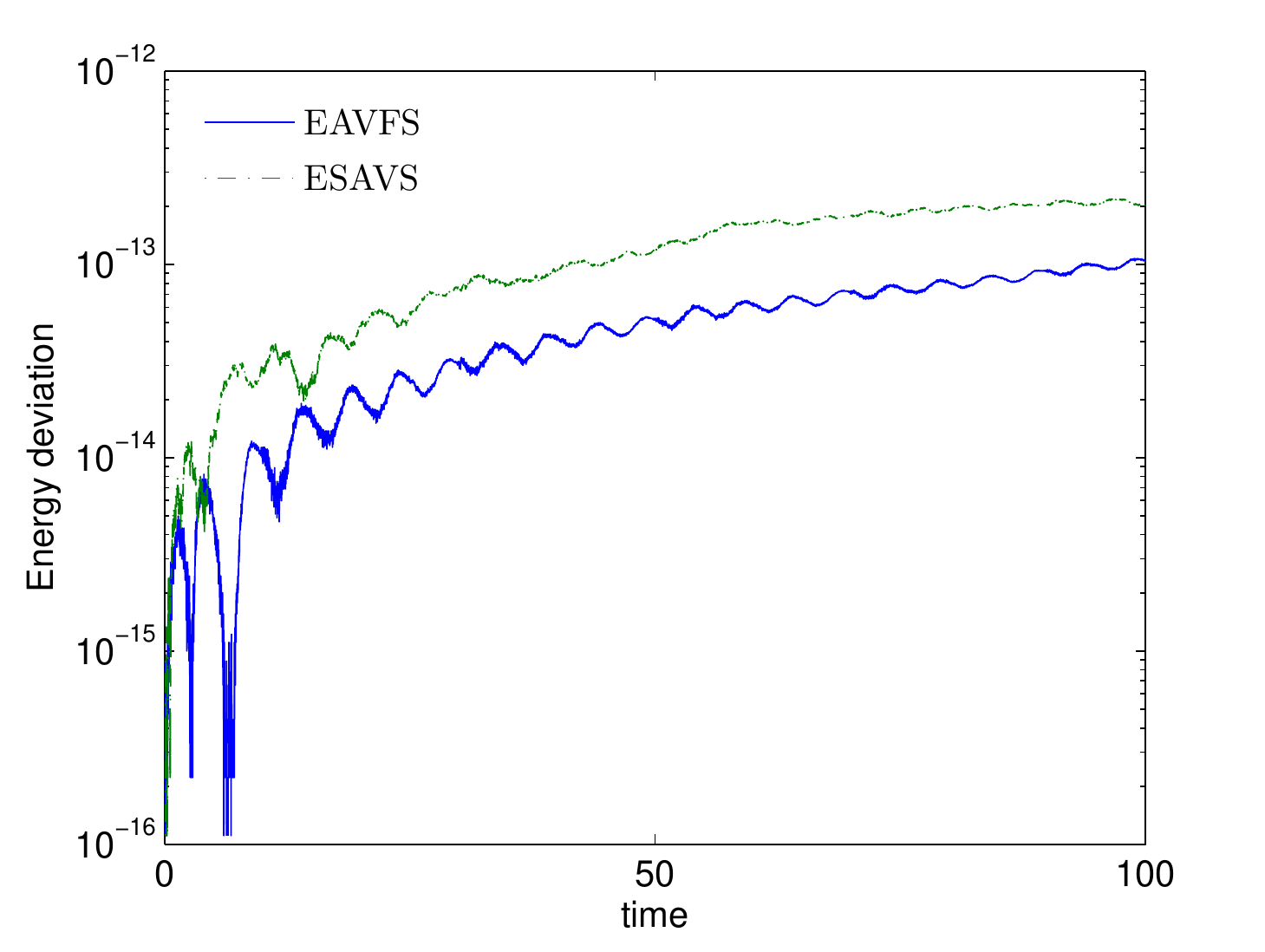}
\end{minipage}
\caption{The energy deviation using the two numerical schemes with time step $\tau=0.01$ and spatial collocation point $N=512$ for
the one dimensional Schr\"odinger equation.}\label{ESAV-scheme:fig:3}
\end{figure}

To further investigate the energy-preservation
of the proposed scheme, we provide the energy errors
using the two numerical schemes for the one dimensional
Schr\"odinger equation over the time interval $t\in[0,100]$ in
Fig. \ref{ESAV-scheme:fig:3}, which shows that all two methods can
exactly preserve the discrete energies.

\begin{figure}[H]
\centering\begin{minipage}[t]{60mm}
\includegraphics[width=60mm]{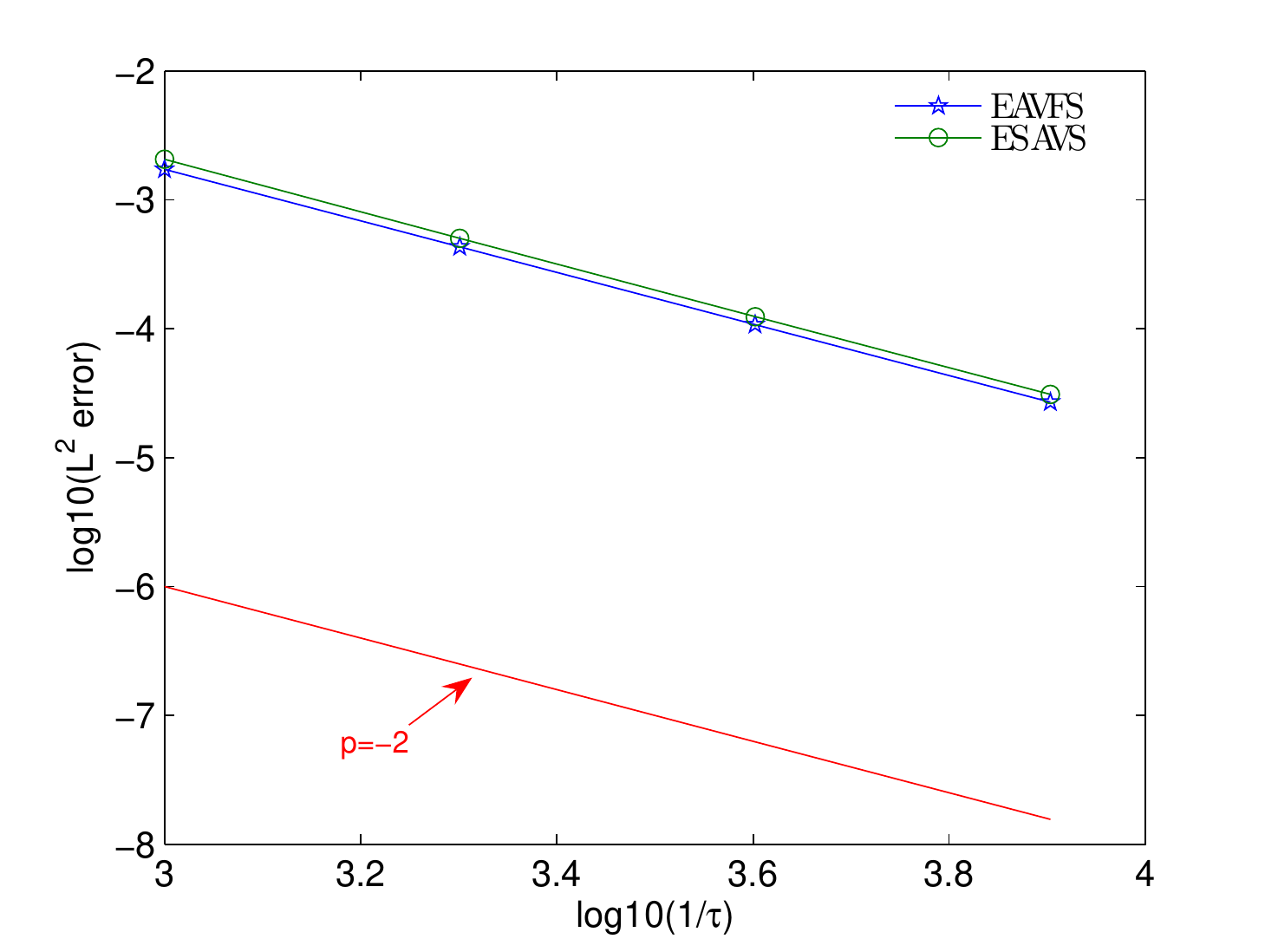}
\end{minipage}
\begin{minipage}[t]{60mm}
\includegraphics[width=60mm]{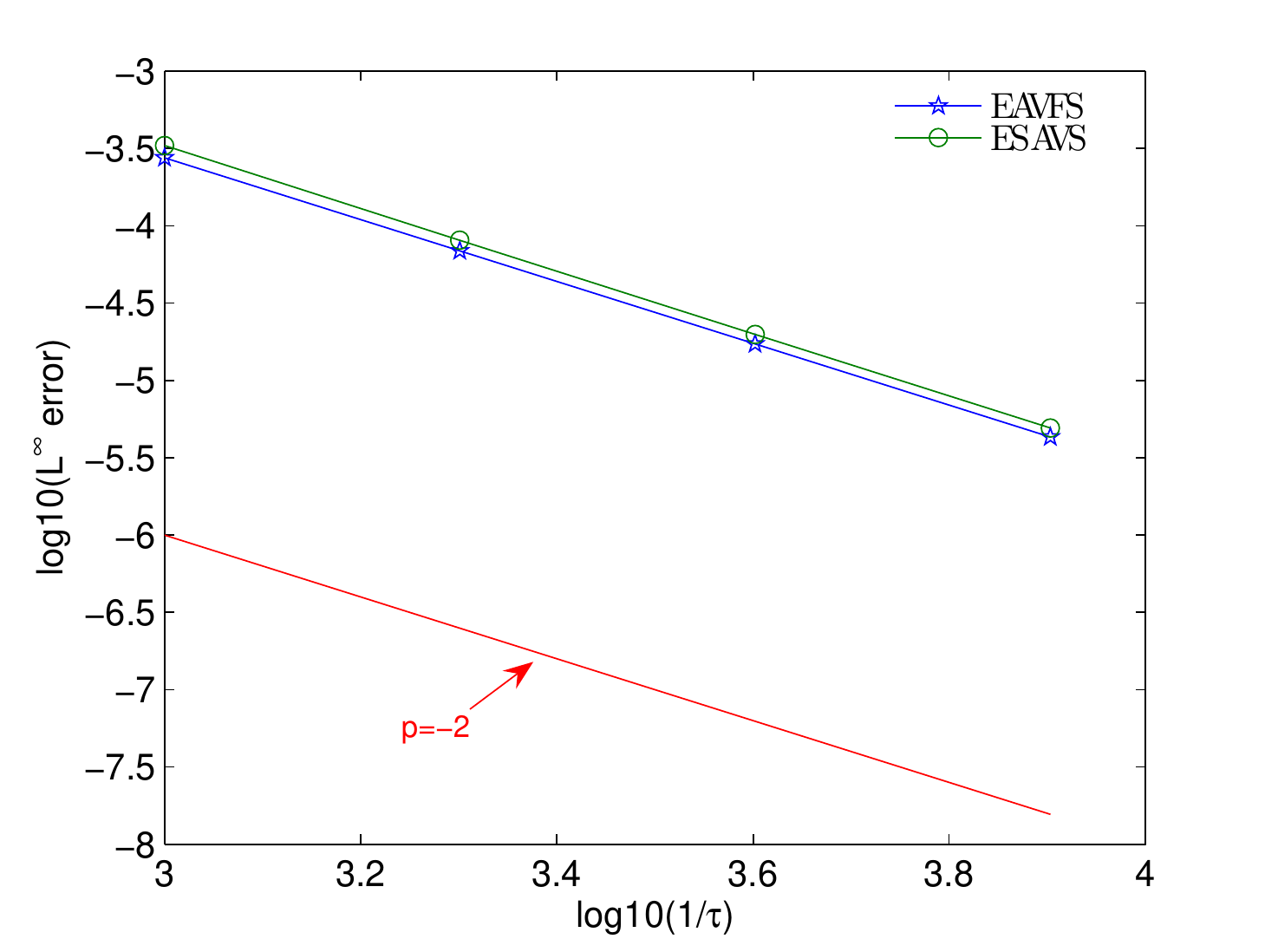}
\end{minipage}
\caption{ Time step refinement tests using the two numerical schemes for the two dimensional Schr\"odinger equation.}\label{ESAV-scheme:fig:4}
\end{figure}

\begin{figure}[H]
\centering\begin{minipage}[t]{60mm}
\includegraphics[width=60mm]{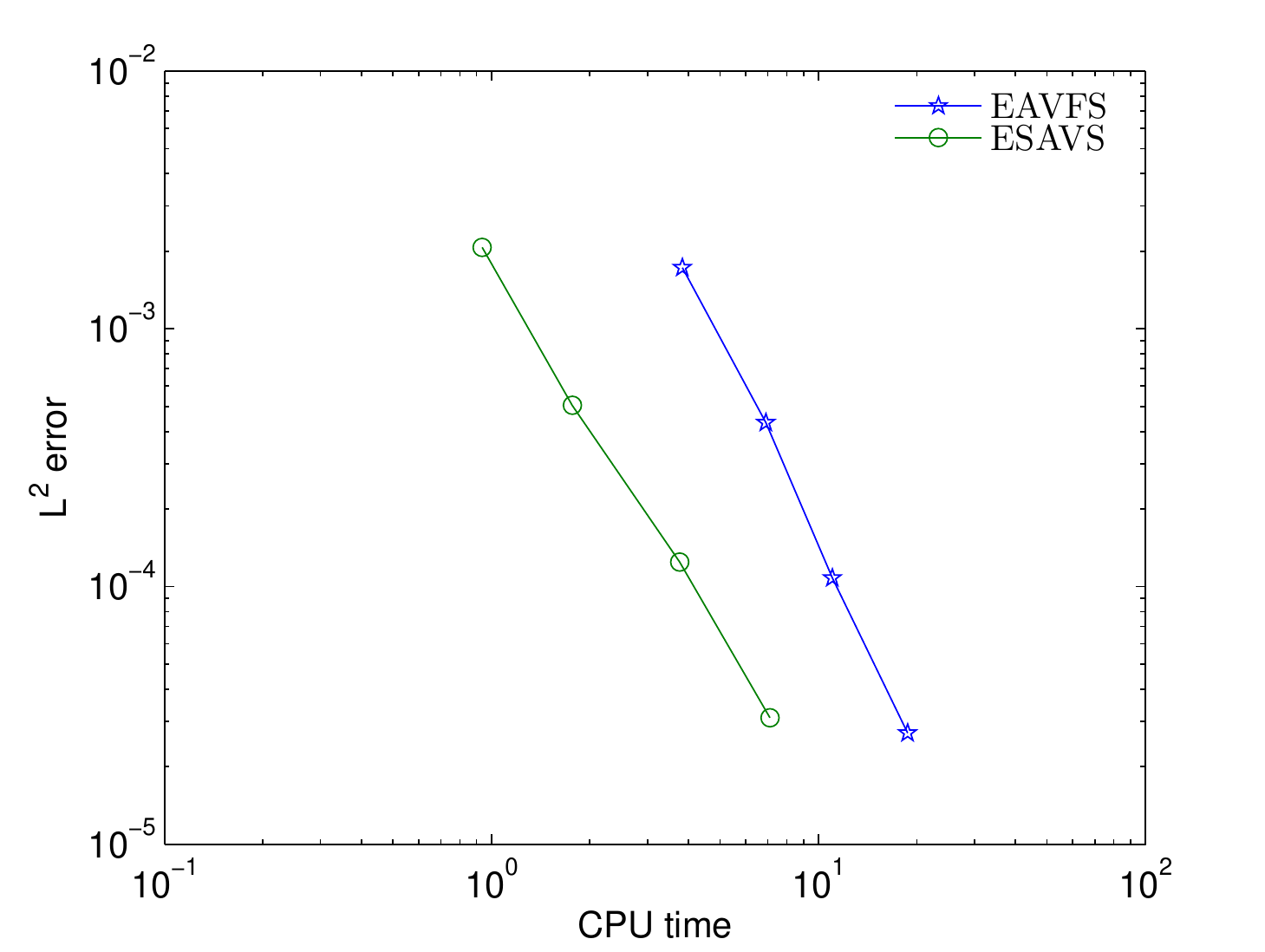}
\end{minipage}
\begin{minipage}[t]{60mm}
\includegraphics[width=60mm]{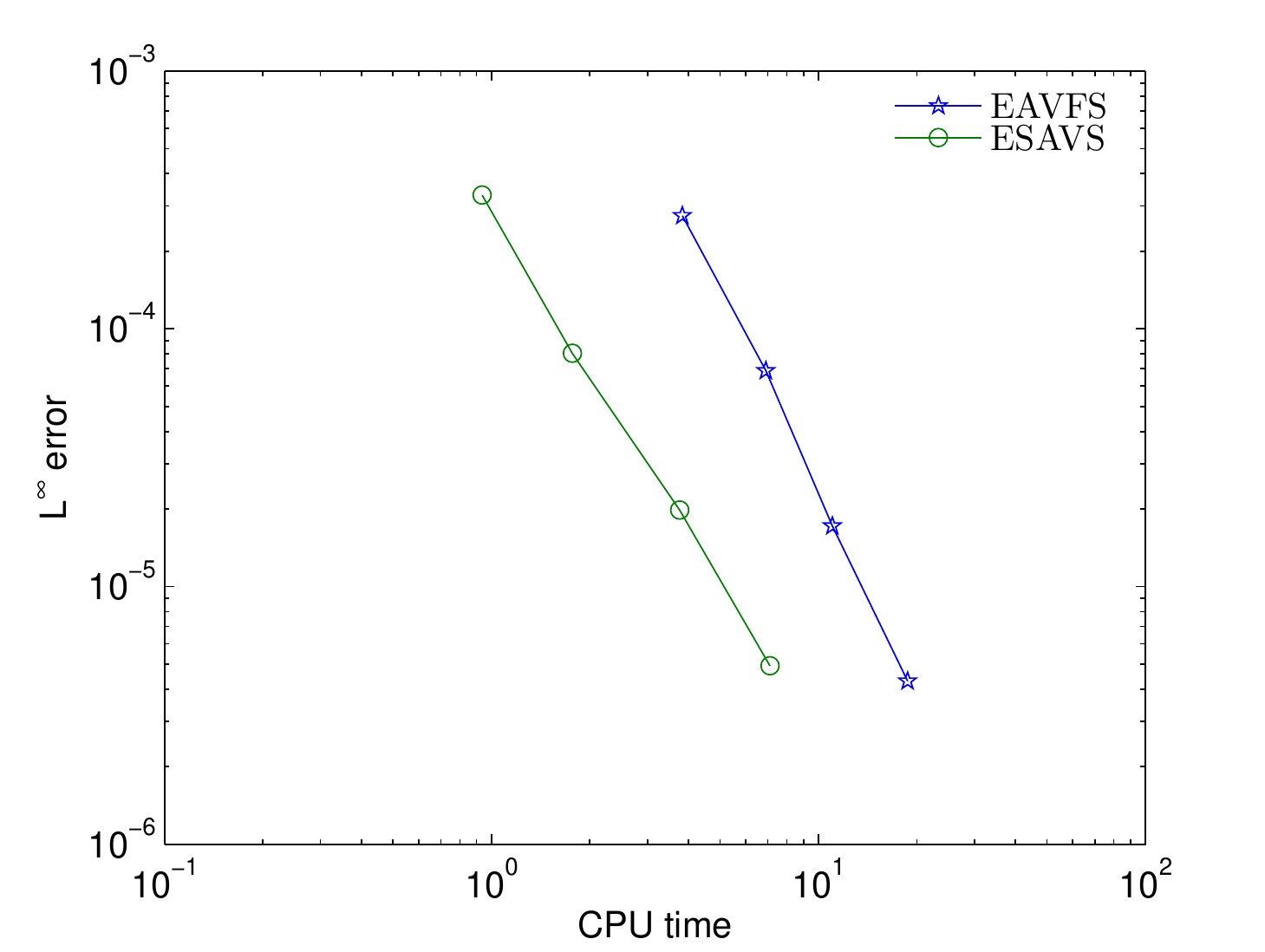}
\end{minipage}
\caption{ The numerical error versus the CPU time using the two numerical schemes for the two dimensional Schr\"odinger equation.}\label{ESAV-scheme:fig:5}
\end{figure}

Next, we apply the proposed scheme to solve the two dimensional nonlinear Schr\"odinger equation which possesses the following analytical solution
\begin{align*}
u(x,y,t)=A\exp(\text{i}(k_1x+k_2y-\omega t)),\ \omega=k_1^2+k_2^2-\beta|A|^2,\ (x,y)\in\Omega.
\end{align*}
We choose the computational domain $\Omega=[0,2\pi]^2$ and take parameters $A=1,\ k_1=k_2=1$, $\beta=-1$ and $C_0=0$. We first test the temporal accuracy of the two numerical schemes by fixing the Fourier node $64\times 64$ such that the spatial discretization errors are negligible. The $L^2$ errors and $L^{\infty}$ errors in numerical solution of $u$ at $t=1$ calculated by using two numerical schemes with various
time steps are shown in Fig. \ref{ESAV-scheme:fig:4}. Also, the global $L^2$ errors and $L^{\infty}$ errors of $u$ versus the CPU time using the two different schemes at $t=1$ are investigated in Fig. \ref{ESAV-scheme:fig:5}. Again, the numerical results indicate that two numerical schemes are second order in time and the proposed scheme is much cheaper than EAVFS.

\begin{figure}[H]
\centering
\begin{minipage}[t]{70mm}
\includegraphics[width=70mm]{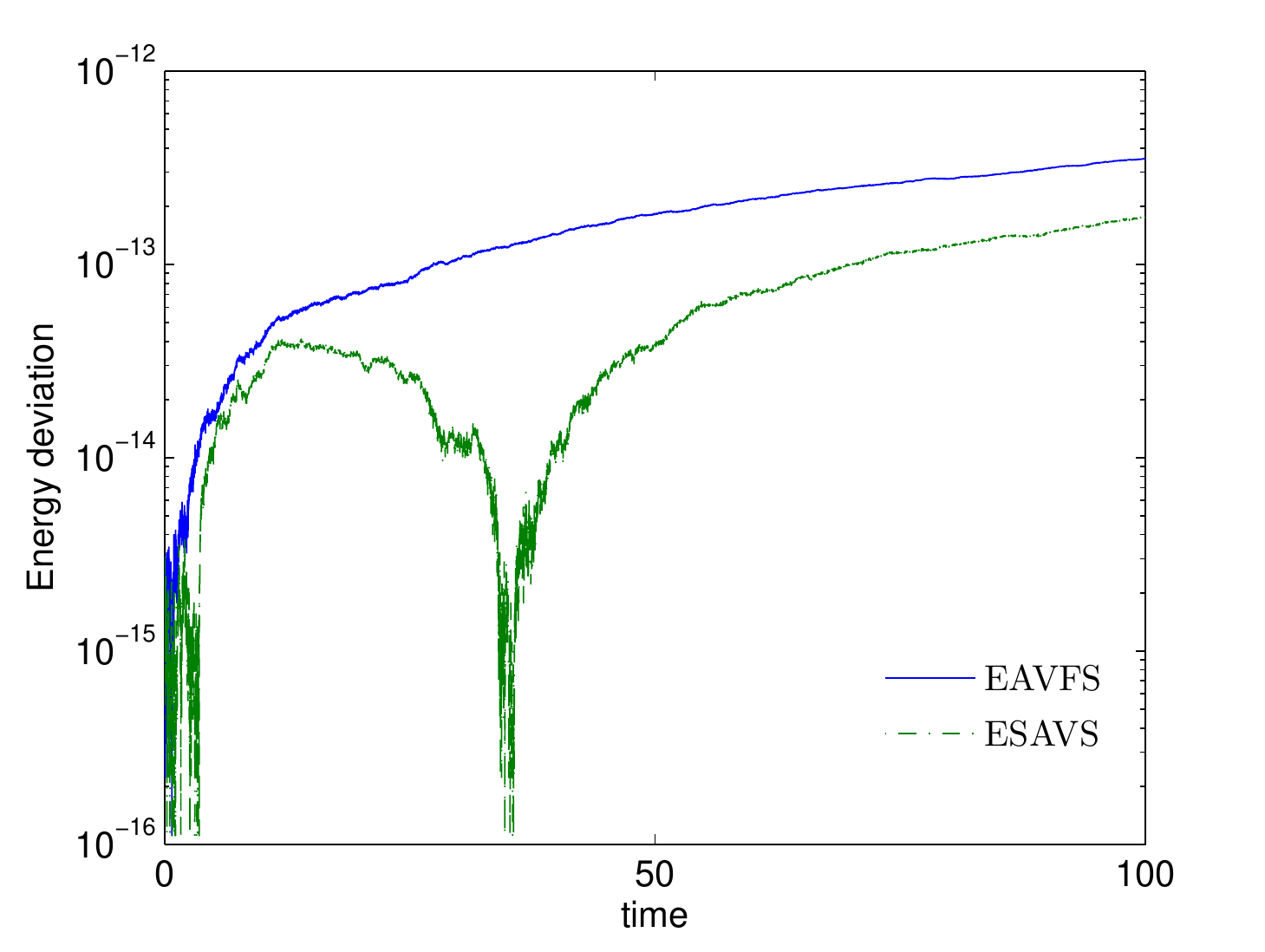}
\end{minipage}
\caption{The energy deviation using the two numerical schemes with time step $\tau=0.01$ and spatial collocation point $32\times 32$ for
the two dimensional Schr\"odinger equation.}\label{ESAV-scheme:fig:6}
\end{figure}

Second, we present the discrete energies for the numerical solutions given by the ESAVS and EAVFS, respectively in Fig. \ref{ESAV-scheme:fig:6}. The numerical
results show that the discrete energies can be exactly preserved, which are consistent with our theoretical result in Theorem \ref{ESAVS-Theorem-5.1}.

\section{Concluding remarks}\label{Sec:NKGEs:6}
In this paper, we design a novel linearly implicit energy-preserving exponential scheme for the nonlinear Klein-Gordon equation and the nonlinear Schr\"odinger equation, respectively. The schemes were developed based on the exponential integrator in combination
with the scalar auxiliary variable (SAV) technique and proved to preserve the discrete energy.  Various numerical examples are carried out to illustrate theoretical analysis. Comparing with the exponential averaged vector filed scheme, the proposed scheme shows remarkable efficiency. Here, we should note that, compared with existing
energy-preserving exponential schemes (e.g., see Refs. \cite{LW16b,SL19jcp}), the proposed method cannot preserve the discrete
Hamiltonian energy. Thus, such trade-offs among methods should
be more carefully investigated. In addition, to the best of our knowledge, the construction of higher order linearly implicit energy-preserving exponential schemes is still not available for conservative systems, which is an interesting topic for future studies.





\section*{Acknowledgments}
The authors would like to express sincere gratitude to the referees for their insightful
comments and suggestions. Chaolong Jiang's work is partially supported by the National Natural Science Foundation
of China (Grant No. 11901513), the Yunnan Provincial Department of Education
Science Research Fund Project (Grant No. 2019J0956) and the Science and Technology
Innovation Team on Applied Mathematics in Universities of Yunnan. Yushun Wang's
work is partially supported by the National Natural Science Foundation of China (Grant
No. 11771213). Wenjun Cai's work is partially supported by the National Natural
Science Foundation of China (Grant No. 11971242) and  the National Key Research and Development Project of China (Grant Nos. 2017YFC0601505, 2017YFC0601406, 2018YFC1504205).


\bibliographystyle{plain}

\end{document}